\documentclass[12pt,reqno]{amsart}

\numberwithin{equation}{section}

\usepackage{amsmath,amsthm,amsfonts,amssymb,amscd,amstext}

\usepackage[pdftex,backref=page]{hyperref}
\usepackage{graphicx,a4wide}
\usepackage{euscript,xcolor}
\usepackage[shortlabels]{enumitem}

\newcommand{\qand}{\quad\text{and}\quad}

\theoremstyle{plain}
\newtheorem{maintheorem}{Theorem}
\newtheorem{maincorollary}[maintheorem]{Corollary}

\newtheorem{theorem}{Theorem}[section]
\newtheorem{proposition}[theorem]{Proposition}
\newtheorem{corollary}[theorem]{Corollary}
\newtheorem{lemma}[theorem]{Lemma}

\theoremstyle{definition}
\newtheorem{remark}[theorem]{Remark}

\newtheorem{conjecture}{Conjecture}

\newcommand{\dem}{\begin{proof}}
  \newcommand{\cqd}{\end{proof}}

\def \epsilon {\varepsilon}

\def \wh {\widehat}
\def \wt {\widetilde}

\def \ov{\overline }

\newcommand{\cC}{\EuScript{C}}
\newcommand{\D}{\EuScript{D}}

\newcommand{\cF}{\EuScript{F}}

\newcommand{\cH}{\EuScript{H}}

\newcommand{\cP}{\EuScript{P}}

\newcommand{\cS}{\EuScript{S}}

\newcommand{\U}{\EuScript{U}}
\newcommand{\cV}{\EuScript{V}}

\newcommand \cU {{\mathcal U}}

\newcommand \vfi {\varphi}

\newcommand{\RR}{{\mathbb R}}

\newcommand{\CC}{{\mathbb C}}
\newcommand{\NN}{{\mathbb N}}

\newcommand{\TT}{{\mathbb T}}
\newcommand{\DD}{{\mathbb D}}
\newcommand{\sS}{{\mathbb S}}

\newcommand{\ZZ}{{\mathbb Z}}

\newcommand{\corr}{\textrm{Cor}}

\newcommand{\interior}{\operatorname{inter}}

\newcommand{\leb}{\operatorname{Leb}}

\newcommand{\dist}{\operatorname{dist}}

\newcommand{\supp}{\operatorname{supp}}

\author{V. Araujo and V. Pinheiro}

\address{V. Araujo, Departamento de Matem\'atica, Universidade Federal da Bahia\\  Av. Ademar de Barros s/n, 40170-110 Salvador, Brazil.}

\email{vitor.d.araujo@ufba.br or vitor.araujo.im.ufba@gmail.com}

\address{V. Pinheiro, Departamento de Matem\'atica,
  Universidade Federal da Bahia, Av. Ademar de Barros s/n,
  40170-110 Salvador, Brazil.}

\email{viltonj@ufba.br}


\date{\today}

\thanks{Work carried out at the Mathematics and Statistics Institute
  of the Federal University of Bahia. V.A. was partially supported by
  CNPq-Brazil (grant 304047/2023-6). V.P. was partially supported by
  IMPA and PRONEX - Dynamical Systems.}

\title[Exponential mixing and the basin problem] {Multidimensional
  non-uniform hyperbolicity, robust exponential mixing and the basin
  problem}

\begin{document}

\begin{abstract}
  We show that the ergodic, topological and geometric basins coincide
  for hyperbolic dominated ergodic $cu$-Gibbs states, solving
  the ``basin problem'' for a wide class of non-uniformly hyperbolic
  systems.
  
  We obtain robust examples of exponential mixing physical measures
  for systems with multidimensional nonuniform hyperbolic dominated
  splitting, without uniformly expanding or contracting subbundles.

  Both results are a consequence of extending the construction of
  Gibbs-Markov-Young structures from partial hyperbolic systems to
  systems with only a dominated splitting, using the existence of an
  ``improved hyperbolic block'', with respect to Pesin's Nonuniform
  Hyperbolic Theory, for hyperbolic dominated measures of smooth maps,
  obtained through hyperbolic times and associated ``coherent
  schedules'' introduced by one of the coauthors. % This also enables us
  % to show that ergodic hyperbolic dominated measures with large
  % entropy have long (un)stable manifolds on positive measure subsets.
  % We also provide examples in the same setting of slower rates of
  % mixing.
\end{abstract}

\maketitle
\tableofcontents

%%%%%%%%%%%%%%%%%%%%%%%%%%%%%%%%%%%%%%%%%%%%%%%%%%%%%%%%%%%%%

\section{Introduction}
\label{sec:introduction}

Dynamical Systems theory is mostly interested in describing the
typical behaviour of orbits as time goes to infinity, and
understanding how this behaviour is modified under small perturbations
of the system. This work concentrates in the study of the former
problem from a probabilistic point of view.  An effective approach is
to describe the average time spent by typical orbits in different
regions of the phase space. According to the Ergodic Theorem (of
Birkhoff), such averages are well defined for almost all points, with
respect to an invariant probability measure. However, frequently the notion of
typical orbit is given in terms  of volume (Lebesgue
measure), which is not always captured by invariant measures. Indeed,
it is a fundamental open problem to understand under which conditions
the behaviour of typical points is well defined, from this statistical
point of view.
% This is always the case if the system preserves Lebesgue measure,
% according to the Ergodic Theorem.

For dissipative systems given by a diffeomorphism $f:M\to M$ on a
phase space $M$, we usually consider the dynamics in the topological
basin of each attracting set, and then restate the question as
follows; see e.g. \cite{BDV2004}.
\begin{description}
\item [Q1]\emph{Is almost every orbit in the basin of attraction
    asymptotic to some orbit contained in the attractor?}
\item[Q2] \emph{Is it generic
    for some natural invariant measure supported in the attractor?}.
\end{description}
An attracting set is a compact invariant subset $A$ of the phase space
$M$ whose \emph{topological basin}
$
  B(A) := \{z\in M:  \omega(x)\subset A \}
$
is a large set --- a neighbourhood of $A$ in our setting (see
e.g.~\cite{Mi85} for other possibilities) --- where
$ \omega(x):=\{y\in M: \exists n_k\nearrow\infty: f^{n_k}(x)\to y\}$
is the set of all accumulation points of the future trajectory of
$x$ (also known as the $\omega$-\emph{limit set} of $x$).

Let $W^s_x$ denote the subset of all points whose trajectory approaches
the future trajectory of the point $x$
\begin{align*}
  W^s_x:=\{ y\in M: \dist(f^ny,f^nx)\to0 \text{  when  } n\nearrow\infty\},
\end{align*}
which, in many cases (e.g. under hyperbolicity assumptions), is a
submanifold of the ambient space.  The \emph{geometric basin} of an
attracting set $A$ is $G(A) := \cup_{x\in A}W^s_x$. We may reformulate
the former question as
\begin{description}
\item[Q1]\emph{does $B(A)=G(A)$ up to zero Lebesgue measure}?
\end{description}
Let us assume that $A$ supports an invariant ergodic probability
measure $\mu$ which is \emph{hyperbolic} (all the Lyapunov exponents
are nonzero) and \emph{physical}, that is, the \emph{ergodic basin}
\begin{align*}
  B(\mu) := \{x\in M: \lim_{n\to+\infty} S_n\vfi(x) /n = \mu(\vfi),
  \forall \vfi\in C^0(M,\RR)\}
\end{align*}
has positive Lebesgue measure in $M$ --- where we denote the ergodic
sum by $S_n\vfi(x):=\sum_{i=0}^{n-1}\vfi(f^ix)$ for any observable
(measurable function) $\vfi:M\to\RR$ and its integral by
$\mu(\vfi):=\int\vfi\,d\mu$. We say that $x\in B(\mu)$ is
$\mu$-generic.  We may now reformulate the latter question as
\begin{description}
\item[Q2]\emph{does $B(\mu) = B(A)$ up to a zero Lebesgue measure set}?
\end{description}
It is well-known that both questions (referred to as ``the basin
problem'') have an affirmative answer in the case of uniformly
hyperbolic (Axiom A) attractors, where the crucial ingredient is the
uniform shadowing property; see e.g.~\cite{Bo75,BR75,Si72,PS89} and
references therein.  On the other hand, not much is known in the
non-uniformly hyperbolic setting: we have positive answers for the
geometric Lorenz-like attracting sets (for which a stable foliation
exists, essentially, by definition~\cite{APPV,AraPac2010}); for
H\'enon-like families from the pioneering work of
Benedicks-Viana~\cite{BeV00} and later
developments~\cite{Tak06,HorMun06} providing strong results on a
nonuniformly hyperbolic setting with no dominated splitting; and for
systems preserving a smooth ergodic measure $\mu$ (where $B(\mu)$ has
full measure as a direct consequence of the ergodic
theorem). Recently, examples of locally dense families of systems with
\emph{historic behavior} (i.e. absence of asymptotic time averages)
for subsets of points with positive volume have been obtained; see
Kiriki-Soma~\cite{KrS15} and together with
Nakano-Vargas~\cite{kiriki2024takens}.

Here we show that \emph{the basin problem always has an affirmative
  answer for hyperbolic dominated $cu$-Gibbs states}, that is,
hyperbolic physical measures admitting a dominated splitting
respecting the hyperbolic decomposition of the Lyapunov exponents,
which are also Sinai-Ruelle-Bowen (SRB) or, equivalently, equilibrium
states with respect to the central-unstable Jacobian.

This is obtained as a consequence of the study of the statistical
properties of physical/SRB measures for non-uniformly hyperbolic
dynamics with a dominated splitting, focusing on the speed of
mixing. For observables (measurable functions) $\vfi,\psi:M\to\RR$,
and an invariant probability measure $\mu$, we consider the
\emph{correlation function}
\begin{align*}
  \corr_{\mu}(\vfi,\psi\circ f^n):=
  |\mu\big( \vfi\cdot \psi\circ f^n\big)  -\mu(\vfi)\mu(\psi)| 
\end{align*}
and recall that $f$ is \emph{mixing} with respect to $\mu$ if
$\corr_\mu(\vfi, \psi\circ f^n) \to 0$ when $n\nearrow\infty$ for any
choice of $\mu$-measurable functions.

In many cases smooth observables satisfy specific rates of decay: in
the uniformly hyperbolic (Axiom A) attractor setting, exponential
mixing holds for H\"older observables with respect to the unique SRB
measure or $u$-Gibbs state~\cite{Bo75,BR75}.  We obtain sufficient
conditions for polynomial and (sub)exponential rates with respect to a
class of $cu$-Gibbs states, which are dominated hyperbolic ergodic
physical measures, using Gibbs-Markov-Young (GMY) structures, as
in~\cite{AlLi15,AlPi10,AlPi08,AlLuPi03}.

These geometric structures were introduced by Young~\cite{Yo98} and
have been applied to study the existence and properties of physical
measures in certain classes of nonuniformly hyperbolic dynamical
systems. GMY structures are known to imply many other statistical
properties beyond the mixing speed, like the Almost Sure Invariance
Principle which then ensures the Central Limit Theorem and the Law of
the Iterated Logarithm~\cite{PhilippStout75}. The speed of mixing is
also strongly related to Large Deviation estimates through the GMY
structure~\cite{AFLV11}.

We extend the construction of these structures from
partially hyperbolic to non-uniformly hyperbolic diffeomorphisms with
a dominated splitting. 

This extension allows us to exhibit examples of robust exponential
mixing for diffeomorphisms without any invariant uniformly hyperbolic
subbundle (expanding or contracting). In our setting the speed of
mixing depends only on the ``tail of hyperbolic times'' along the
central unstable direction. We note that
Melbourne-Varandas~\cite{melvar16} showed that exponential contraction
(and expansion) along the stable (and unstable) direction, at the
returns of a generalized horseshoe on a well-chosen subset of the
ambient space, is enough to build GMY structures.

Here, we do not need to assume any condition on the speed of
convergence of non-uniform contraction along the center-stable
direction to obtain specific rates of mixing, since we obtain
\emph{uniformly long stable leaves with uniform contraction Lebesgue
  almost everywhere inside certain cylinders on the ambient space} ---
providing the ``generalized horseshoe with infinitely many returns in
variable times'' as in Young~\cite{Yo98} --- which, in turn, enables a
solution to the basin problem.

We use the existence of an ``improved hyperbolic block'' (akin to the
hyperbolic blocks of the Nonuniform Hyperbolic Theory of
Pesin~\cite{BarPes2007}) for hyperbolic dominated measures of smooth
maps, obtained through hyperbolic times and associated ``coherent
schedules'', as a sharp tool to prove our results. % This also enables
% us to show that ergodic hyperbolic dominated measures with large
% entropy have long (un)stable manifolds on positive measure subsets,
% akin to the recent results from Luo-Yang~\cite{LuoYang24} with the
% extra domination assumption but without the need of $C^\infty$
% regularity.

% \hrulefill

% Henon BC attractor: example of exponential mixing without domination

% History Rob. Mix: Hyp. Attractors, part. hyperbolic diffeos, flows --
% Lorenz-like, Anosov flows

% enunciar teoremas gerais com dominacao no suporte de medida SRB para
% obter GMY e velocidades de mistura

% observar que blocos coerentes são mais fortes do que blocos
% hiperbolicos de Pesin ``classicos''

% \hrulefill

\section{Statement of results}
\label{sec:statement-results}

Let $M$ be a compact finite dimensional Riemannian manifold with an
induced distance $d$ and volume form $\leb$. If $M$ has a boundary,
then we assume that all the maps $f:M\to M$ to be considered send the
boundary in the interior $f(\partial M)\subset M\setminus\partial M$,
in what follows.

% Given $x\in M$ we denote the \emph{omega-limit} set of accumulation
% points of the future orbit of $x$ as
% $ \omega(x):=\{y\in M: \exists n_k\nearrow\infty: f^{n_k}(x)\to y\}.
% $ Given a compact invariant subset $ A\subset M$ we denote by
% $\M( A,f)$ the family of all $f$-invariant probabilty measures
% supported on $ A$ and by $\M(f)$ the set of all $f$-invariant
% probability measures.

% We also denote by $\M(x)$ the set of $f$-invariant probability
% measures which are weak$^*$ accumulation points of the empirical
% measures $\mu_n(x):=n^{-1}\sum_{i=0}^{n-1}\delta_{f^ix}$. That is,
% $\nu\in\M(x)$ means that we can find a sequence of positive integers
% $n_k\nearrow\infty$ so that $ \int\vfi\,d\mu_n(x) \to \int\vfi\,d\nu $
% for any continuous function $\vfi:M\to\RR$.

Let $f:M\circlearrowleft$ be a diffeomorphism and $ A$ a compact
$f$-invariant subset. We say that $ A$ has a dominated splitting
if there exists an $Df$-invariant splitting
$ T_A M=E^{cs}_A\oplus E^{cu}_A$ and constants
$0<\lambda<1$, $c>0$ such that for all $n\ge1$ and $x\in A$
\begin{align*}
  \|Df^n\mid E^{cs}_x\|\cdot\|(Df^n\mid E^{cu}_{f^nx})^{-1}\|\le
  c\lambda^n.
\end{align*}

\subsection{Non-uniform expansion and/or contraction, Gibbs states and
  physical measures}
\label{sec:non-uniform-expans}

The following notions imply non-negative Lyapunov exponents and have
been used to obtain physical measures and study their statistical
properties since~\cite{BoV00,ABV00}.

\subsubsection{Non-uniform hyperbolicity}
\label{sec:non-uniform-expans-1}

For any function $\vfi:M\to\RR$ and map $g:M\circlearrowleft$ we write
$S^g_n\vfi$ for the ergodic sum $\sum_{i=0}^{n-1}\vfi\circ g^i$.
We set $\phi_k^{cu}(x):=\log\|(Df^k\mid E^{cu}_x)^{-1}\|$ and
$\phi_k^{cs}(x):=\log\| Df^k\mid E^{cs}_x \|$ for each $k\ge1$ in what
follows and write $\phi^*=\phi_1^*$ for $*=cs,cu$.

We say that the center-unstable subbundle $E^{cu}$ is
\emph{non-uniformly expanding} (with respect to $\leb$) if we can find
$c_u>0$ and a subset $H_u$ with $\leb(H_u)>0$  so that
\begin{align}\label{eq:NUE0}
  \limsup\nolimits_{n\to\infty} S_n\phi^{cu}(x)/n < -c_u,
  \quad\text{for } x\in H_u.
\end{align}
We say that the center-stable bundle $E^{cs}$ is \emph{non-uniformly
  contracting} (with respect to $\leb$) if we can find $c_s>0$ and a
subset $H_s$ with $\leb(H_s)>0$ so that
\begin{align}\label{eq:NUC}
  \limsup\nolimits_{n\to\infty}S_n\phi^{cs}(x)/n<-c_s, \quad\text{for
  $x\in H_s$}.
\end{align}
We say that a diffeomorphism $f$ with a globally defined dominated
splitting, whose bundles are both non-uniformly expanding and
contracting (with respect to $\leb$) on the same $\leb$-positive
subset $H := H_u\cap H_s$, is \emph{non-uniformly hyperbolic}.

\subsubsection{Hyperbolic and dominated invariant probability measures}
\label{sec:invariant-sets-with}

An $f$-in\-variant probability measure $\mu$ is
\emph{hyperbolic} if the Lyapunov exponents provided by Oseledets'
Multiplicative Ergodic Theorem $\mu$-a.e. are all non-zero. We say
that $\mu$ is \emph{hyperbolic and dominated} if its support
$\supp\mu$ admits a dominated splitting $E^{cs}\oplus E^{cu}$ which
separates the hyperbolic Oseledets subspaces in the following sense:
for $\mu$-a.e. $x$
\begin{align}\label{eq:middlexp}
  \lambda_{cs}^+:=\lim_{n\to+\infty}\log\|Df^n\mid E^{cs}_x\|^{1/n}<0
  \;\&\;
  \lambda_{cu}^-:=\lim_{n\to+\infty}\log\|(Df^n\mid E^{cu}_x)^{-1}\|^{1/n}<0.
\end{align}

\subsubsection{Attracting sets}
\label{sec:attracting-sets}

We say that an invariant subset $A$ is \emph{attracting} if there
exists a \emph{open trapping neighborhood} $U$ of $A$ so that
$\ov{f^k(U)}\subset U$ for some $k\ge1$ and
$A=\cap_{n\ge1}\ov{f^n(U)}$. If additionally $A$ admits a dense
forward trajectory, that is, if we can find $x\in A$ so that
$\omega(x)=A$, then $A$ is an \emph{attractor}.

If $A$ admits a dominated splitting, then we can extend the splitting
continuously to a small neighborhood $U$ of $A$. We may assume without
loss of generality that $U$ is a trapping neighborhood.

We say that an attracting set $A$ with a dominated splitting is
non-uniformly hyperbolic (with respect to $\leb$) if the extended
bundles satisfy both~\eqref{eq:NUE0} and~\eqref{eq:NUC} on the same
$\leb$-positive measure subset $H\subset H_s\cap H_u\subset U$.

\subsubsection{Gibbs states}
\label{sec:gibbs-states}

We say that an $f$-invariant probability measure $\mu$ supported on a
compact invariant subset $A$ with dominated splitting is a
\emph{$cu$-Gibbs state} if
\begin{enumerate}[(i)]
\item $\mu$ satisfies the \emph{Entropy Formula}: if $h_\mu(f)$ is the
  Kolmogorov-Sinai entropy of the measure preserving system
  $(M,f,\mu)$ and $J^{cu}:=\log|\det(Df\mid E^{cu})|$ is the logarithm
  of central-unstable Jacobian, then
$
    h_\mu(f)=\int J^{cu}\,d\mu;
 $
\item all Lyapunov exponents along $E^{cu}$ are positive $\mu$-almost
  everywhere, that is, for
  $\mu$-a.e. $x$ we have
  % , we can find $n>1$ such that
  $\lim_{n\nearrow\infty}\log\|(Df^n\mid E^{cu})^{-1}\|^{1/n}<0$.
% $\mu(\phi^{cu}_n)=\int \log\|(Df^n\mid E^{cu})^{-1}\| \,d\mu <0$.
  \end{enumerate}

  \begin{remark}[Hyperbolic dominated Gibbs states and non-uniform
    hyperbolicity]
    We recall that if $\mu$ is a hyperbolic $cu$-Gibbs state with
    dominated splitting, then some power $g:=f^N$ is non-uniformly
    hyperbolic, that is, both conditions~\eqref{eq:NUE0}
    and~\eqref{eq:NUC} hold for $\mu$-a.e. $x$ with respect to $g$,
    and so $\leb(H)>0$; see Theorem~\ref{mthm:Gibbsmix} and
    Subsection~\ref{sec:hyperb-dominat-measu-1} and cf.~\cite{ABV00}.
  \end{remark}

 \subsection{Ergodic and geometric basin coincide Lebesgue
   modulo zero}
\label{sec:ergodic-basin-geomet}

The following extends the positive answer to the basin problem from
uniformly hyperbolic (Axioma A) $C^2$ diffeomorphisms to a much wider
class of smooth nonuniformly hyperbolic systems.

\begin{maintheorem}
  \label{mthm:basin}
  Let $g:M\circlearrowleft$ be a $C^{1+\eta}$ diffeomorphism, for some
  $\eta\in(0,1]$, with a dominated splitting
  $T_AM=E^{cs}_A\oplus E^{cu}_A$ over an attracting set $A$ on a
  trapping neighborhood $U\subset M$, and an ergodic hyperboliic
  dominated $cu$-Gibbs state $\mu$ for $g$ with $\supp(\mu)\subset
  A$. Then modulo zero volume subsets we have
  \begin{align*}
    %\label{eq:basin}
     G(\supp \mu) = B(\mu).
  \end{align*}
  If $A$ is an attractor (i.e.,
  transitive), then $\supp\mu=A$ and we obtain
  \begin{align*}
    B(A)=G(A)=B(\mu)
  \end{align*}
    modulo zero volume subsets.
  \end{maintheorem}
  
The proof is a scholium of the study of statistical properties of such
invariant measures whose results we present in what follows.

\begin{remark}[wild attractors]
  \label{rmk:wildatt}
  This shows that the class of attractors in the statement of
  Theorem~\ref{mthm:basin} \emph{are not wild}.  We recall that a
  \emph{wild attractor} admits a cycle of subsets
  $A=A_0\cup\dots\cup A_{s-1}$ for some $s\ge1$ so that
  $f(A_i)=A_{(i+1)\bmod s}, i\ge0$ and $f\mid_A$ is transitive; but
  there exists a (Cantor) subset $\Lambda\subset A$ so that
  $\omega(x)=\Lambda$ for $\leb$-a.e. $x\in A$; see
  e.g.~\cite{Mi85,BKNvS96}.
\end{remark}

\subsection{Average expansion times and mixing for hyperbolic
  dominated Gibbs state}
\label{sec:exponent-mixing-phys}

% We say that the center-stable bundle $E^{cs}$ is \emph{weakly
%   non-uniformly contracting} on $H$ if we can find $c_s>0$ and for
% each $x\in H$ a strictly increasing subsequence $n_k=n_k(x)$ so that
% for some $b>0$
% \begin{align}\label{eq:wNUC}
%   \limsup (n_{k+1}/n_k)> b
%   \qand
%   S_{n_k}\phi^{cs}(x)<-c_s n_k, \quad k\ge1.
% \end{align}
% and define the \emph{average contracting time} function
% \begin{align*}
%   c_H(x)
%   =
%   \min\{N\ge1: S_n\phi^{cs}(x) < -nc_s\}
% \end{align*}
Given any embedded disk $\Sigma$ in $M$ we denote by $\leb_\Sigma$ the
induced volume form on $\Sigma$. From the existence of the dominated
splitting, for small $a>0$ we find center unstable and stable cones
\begin{align}
  \label{eq:cucone}
  C^{cu}_a(x)
  &=
    \{ v=v^s+v^c : v^s\in E^{cs}_x, v^c\in E^{cu}_x, x\in
    M, \|v^s\|\le a\|v^c\|\}, \qand
  \\
  C^{cs}_a(x)
  &=
    \{ v=v^s+v^c : v^s\in E^{cs}_x, v^c\in E^{cu}_x, x\in
    M, \|v^c\|\le a\|v^s\|\}, \nonumber
\end{align}
which are invariant in the following sense
\begin{align}\label{eq:coneinv}
  Df(x)\cdot C^{cu}_a(x)\subset C^{cu}_a(f(x))
  \qand
  Df\cdot C^{cs}_a(x)\supset C^{cs}_a(f(x)),
\end{align}
for all $x\in U$.  We say that an embedded $C^1$ disk $\Sigma$ is a
\emph{$cu$-disk} if $T_x\Sigma\subset C^{cu}_a(x)$ for all
$x\in\Sigma$ (and, analogously, a \emph{$cs$-disk} if
$T_x\Sigma\subset C^{cs}_a(x)$ for all $x\in \Sigma$).

Putting together the main results of this text and other known
standard results of non-uniform hyperbolic dynamics, we obtain the
following.

\begin{maintheorem}
  \label{mthm:Gibbsmix}
  Let $f:M\circlearrowleft$ be a $C^{1+\eta}$ diffeomorphism, for some
  $\eta\in(0,1]$, with a dominated splitting
  $T_AM=E^{cs}_A\oplus E^{cu}_A$ over an attracting set $A$ on a
  trapping neighborhood $U\subset M$, admitting an ergodic hyperbolic
  dominated $cu$-Gibbs state $\mu$.
  Then
  \begin{enumerate}[(A)]
  \item there exists $N\ge1$ such that $g:=f^N$ is non-uniformly
    expanding along the center-unstable direction and non-uniformly
    contracting along the center-stable direction with respect to
    Lebesgue measure.
  \end{enumerate}
  Let $H\subset M$ be the subset of points $x\in M$ where non-uniform
  hyperbolicity holds and define the \emph{expansion time function}
  (which is finite for the points $x\in H$)
\begin{align}\label{eq:cutail}
  h(x)
  =
  h^{cu}(x)
  =
  \min\left\{
  N\ge1: S^g_n\phi_N^{cu}(x) < -nc_u/2, \quad \forall\, n\geq N\right\}.
\end{align}
Then we can find and integer $q\ge1$ so that $g^q$ has $1\le p\le q$
invariant mixing probability measures $\nu_1,\dots,\nu_p$ so that
$f_*\nu_i=\nu_{i+1}$ for $i=1,\dots, p-1$; $f_*\nu_p=\nu_1$, and
$\mu=\frac1p\sum_{i=1}^p \nu_i$. In addition, for each $1\le i\le p$
  \begin{enumerate}[(A),resume]
  \item if, moreover, for some $cu$-disk $\gamma\subset A$
    admitting a full $\leb_\gamma$-measure subset of $\mu$-generic
    points\footnote{It follows from the construction of GMY structure
      that these disks always exist; see
      Subsection~\ref{sec:gibbsm-struct} and
      Remark~\ref{rmk:fulldiskH}.}, the expansion time function $h$
    for the dynamics of $g$ satisfies
    \begin{enumerate}[(1)]
    \item $\leb_\gamma\{h\ge n\}\le Cn^{-\alpha}$ for some $C>0$ and
      $\alpha>1$, then $(g^q,\nu_i)$ mixes polynomially, i.e., for
      all $\eta$-H\"older observables $\vfi,\psi:M\to\RR$ there is
      $C'>0$ so that
      $\corr_{\nu_i}(\vfi,\psi\circ g^{qn})\le C'n^{-\alpha+1}$ for
      all $n\ge1$;
    \item $\leb_\gamma\{h\ge n\}\le C e^{-c n^\alpha}$ for some
      $C,c>0$ and $0<\alpha\le1$, then $(g^q,\nu_i)$ mixes
      (sub)exponentially, i.e., there exists $c'>0$ such that
      $\eta$-H\"older observables $\vfi,\psi:M\to\RR$ admit $C'>0$ for
      which
      $\corr_{\nu_i}(\vfi,\psi\circ g^{qn})\le C'e^{ - c' n^\alpha}$
      for all $n\ge1$.
    \end{enumerate}
  \end{enumerate}
\end{maintheorem}

\begin{remark}
  \label{rmk:nocontrolcs}
  There is no need of control hiperbolicity along the center-stable
  direction.
\end{remark}

\subsubsection{Robust non-uniformly hyperbolic exponentially mixing
  class}
\label{sec:robust-non-uniformly}

We recall that $f$ is \emph{topologically mixing over an invariant
  subset} $A$ if for each pair of nonempty open subsets $U,V$ so that
$U\cap A\neq\emptyset\neq V\cap A$ there exists $N>1$ such that
$V\cap f^n U \neq0$ for all $n>N$.

\begin{maincorollary}
  \label{mcor:GMYexp}
  In the same setting of Theorem~\ref{mthm:NUHypGMY},  if we
  additionally assume that:
  \begin{itemize}
  \item $f$ is topologically mixing over $A$; and 
  \item admits a $cu$-disk $\gamma$ contained in $A$, such that
    $\gamma$ contains a full $\leb_\gamma$-measure subset of
    non-uniformly hyperbolic points, satisfying
    $\leb_\Sigma(h>n)\le C e^{-n\zeta}$ for some $C,\zeta>0$ and all
    $n>1$.
  \end{itemize}
  Then there exists $\omega>0$ so that, for any $\eta$-H\"older
  observables $\psi_1,\psi_2$, we can find $C'>0$ so that
  $C_\mu(\psi_1,\psi_2)\le C'e^{-n\omega}$ for all $n\ge1$.
\end{maincorollary}

% see e.g.~\cite{MCY2017} for conditions on the existence of
% physical/SRB measure on a slightly more general setting and
% also~\cite{ABV00,ADLP}.

%The following is the motivating example for the results in this work.

%\begin{remark}[Robust exponencial mixing]
The robust topologically mixing class of $C^2$ diffeomorphisms on the
$n$-torus ($n\ge4$ with $A=M$) from Tahzibi~\cite{tah04}, described in
the following Section~\ref{sec:auxiliary-results} (see
Theorem~\ref{thm:NUEnotPH}), together with Corollary~\ref{mcor:GMYexp}
provide the existence of \emph{robust non-uniformly hyperbolic
  exponentially mixing diffeomorphisms} (from
Proposition~\ref{pr:robustexpmix}), \emph{without any uniformly
  contracting or expanding invariant subbundle}.
%\end{remark}

\subsection{GMY structure for hyperbolic dominated $cu$-Gibss states}
\label{sec:gibbs-markov-young}

We show that all $cu$-Gibbs states which are hyperbolic and dominated
must have a GMY structure with integrable return times, which enables
us to study mixing rates for these types of invariant probability
measures, as in Theorem~\ref{mthm:Gibbsmix}.

% Assuming that
% non-uniform expansion~\eqref{eq:NUE0} together with
% non-uniform contraction~\eqref{eq:NUC} holds on a full Lebesgue
% measure subset of  a $cu$-disk $\Sigma$ which is also subset of
% Birkhoff generic points of a physical/SRB ergodic probability measure
% $\mu$, we are able to show the existence of a very useful structure
% enabling us to study statistical properties of $\mu$.

% with a full $\leb_\Sigma$ subset of points satisfying simultaneously
% non-uniform expansion and non-uniform contraction with uniformly sized
% stable leaves,
%leading to the following conclusion.

\begin{maintheorem}
  \label{mthm:NUHypGMY}
  Let $f:M\circlearrowleft$ be a $C^{1+\eta}$ diffeomorphism, for some
  $\eta\in(0,1]$, with a dominated splitting
  $T_AM=E^{cs}_A\oplus E^{cu}_A$ over an attracting set $A$ on a
  trapping neighborhood $U\subset M$, and an ergodic hyperbolic
  dominated $cu$-Gibbs state $\mu$ for $f$.
  % , for which there is an $cu$-disk
  % $\Sigma$ that is non-uniformly expanding along the $E^{cu}$
  % direction and non-uniformly contracting along the $E^{cs}$ direction
  % on a subset $H\subset\Sigma$ with full $\leb_\Sigma$-measure of
  % $\mu$-generic points.

  Then, for some $k\ge1$, $g=f^k$ admits a GMY structure
  $\Lambda\subset A$ for $\mu$ with integrable return times.
\end{maintheorem}

For the detailed definition of a GMY structure, see
Section~\ref{sec:existence-physic-mea-1}.
These geometric structures were introduced by Young~\cite{Yo98} and
have been applied to study the existence and properties of physical
measures in certain classes of nonuniformly hyperbolic dynamical
systems.

Theorem~\ref{mthm:NUHypGMY} is an extension of~\cite[Corollary
7.28]{Alves2020b} from a partially hyperbolic non-uniformly expanding
setting to the setting of dominated splitting with non-uniform
hyperbolicity \emph{with the extra assumption of existence of a
  hyperbolic $cu$-Gibbs state}.

\subsubsection{Existence of physical measures and GMY structure}
\label{sec:existence-physic-mea}

The non-uniform hyperbolic assumption on $A$, as in~\eqref{eq:NUE0}
and~\eqref{eq:NUC} with $\leb(H)>0$, does not ensure that all ergodic
$cu$-states are hyperbolic (or physical measures); see
Remark~\ref{rmk:nonsuperfl}.

The existence of ergodic hyperbolic $cu$-Gibbs states in our setting
can be ensured under an extra assumption. We say that $f$ is
\emph{mostly contracting} along the center-stable subbundle if
\begin{align}
  \label{eq:mostcont}
  \limsup\nolimits_{n\nearrow\infty}\log\|Df^n\mid E^{cs}_x\|^{1/n}<0
\end{align}
for \emph{a positive Lebesgue measure set of points $x$ in every $cu$-disk
inside} $U$. 

\begin{theorem}{\cite[Theorem C]{vasquez2006}}
  \label{thm:Vasquez}
  Let $f:M\circlearrowleft$ be a $C^{1+\eta}$ diffeomorphism, for some
  $\eta\in(0,1]$, with a dominated splitting
  $T_AM=E^{cs}_A\oplus E^{cu}_A$ over an attracting set $A$ on a
  trapping region $U\subset M$, which is nonuniformly expanding along
  $E^{cu}$ and mostly contracting along $E^{cs}$.
  Then $f$ admits finitely many ergodic physical/SRB
  measures $\mu_1,\dots,\mu_k$ which are $cu$-Gibbs states and whose
  basis cover $\leb$-a.e point of $H$, that is: for each $i=1,\dots,k$
  the ergodic basin of $\mu_i$ 
  % \begin{align*}
  %   B(\mu_i)=\left\{x\in
  %   M:\frac1nS_n\vfi(x)\xrightarrow[n\to\infty]{}\int\vfi\,d\mu,\,
  %   \forall\vfi\in C^0(M,\RR)\right\}
  % \end{align*}
  has positive volume $\leb(B(\mu_i))>0$, and
  $\leb\big(H\setminus(B(\mu_1)\cup\dots\cup B(\mu_k))\big)=0$.
\end{theorem}

We obtain the following improvement of the results from
Alves-Bonatti-Viana~\cite{ABV00} and Vasquez~\cite{vasquez2006}.  We
say that an attracting set $A$ is \emph{weakly dissipative} if it
admits a trapping neighborhood $U$ so that
$J(x):=\log |\det Df_x|\le0$ for all $x\in U$.

\begin{maincorollary}
  \label{mcor:abv}
  Every non-uniformly hyperbolic \emph{weakly dissipative} attracting
  set of a $C^{1+}$-diffeomorphism $f$ with \emph{one-dimensional
    center-stable bundle} satisfies the same conclusion of
  Theorem~\ref{thm:Vasquez}. Moreover, each physical/SRB measure
  $\mu_i$ admits a GMY structure with integrable return times. In
  addition, if $\leb(U\setminus H)=0$, then we get
  \begin{align*}
    B(A)= B(\mu_1)\cup\ldots\cup B(\mu_k)=G(A), \quad \leb-\bmod0.
  \end{align*}
\end{maincorollary}

\begin{remark}
  \label{rmk:conformal}
  We may replace the assumption of one-dimensional center-stable
  bundle by a conformal center-stable bundle with any finite dimension
  and keep the conclusion of Corollary~\ref{mcor:abv}, that is, we may
  assume that $Df(x)v=a(x)\cdot v$ for each $v\in E^{cs}_x$, $x\in M$
  where $a:\Lambda\to\RR$ is H\"older-continuous. Without
  conformality, see Conjecture~\ref{conj:wNUHyphys} in the following
  Subsection~\ref{sec:comments-conjectures}.
\end{remark}

\subsection{Consequences for hyperbolic dominated measures}
\label{sec:conseq-hyperb-domina}

We now consider ergodic hyperbolic dominated invariant probability
measures which are not necessarily $cu$-Gibbs states. The statement of
the next theorem assumes the usual non-uniform hyperbolic condition
from Pesin's Theory plus domination, and provides a ``hyperbolic
coherent block'' with positive measure and strong uniformly
hyperbolic features.

To present the next result, we say that an embedded disk
$\gamma\subset M$ is a (local) \emph{unstable manifold}, or an
\emph{unstable disk}, if $d(f^{-n}(x),f^{-n}(y))$ tends to zero
exponentially fast as $n\nearrow\infty$, for every
$x,y\in\gamma$. Analogously, $\gamma$ is a (local) \emph{stable
  manifold}, or a \emph{stable disk}\footnote{Cf. the definition of
  $cu$-disk and $cs$-disk before the statement of
  Theorem~\ref{mthm:NUHypGMY}.}, if $d(f^n(x),f^n(y))\to0$
exponentially fast as $n\nearrow\infty$, for every $x,y\in\gamma$. We
say that $\gamma$ has inner radius larger than $\delta>0$ around $x$,
if there exists a closed $\delta$-neighborhood $T^\delta_x$ of the
origin in $T_x\gamma$ and an immersion $i:T^\delta_x\to\gamma$ so that
the intrinsic distance between $i(0)$ and $i(p)$ within $\gamma$, for
any $p\in\partial T^\delta_x$, is at least $\delta$.

\begin{maintheorem}[Long (un)stable leaves with positive frequency]
  \label{mthm:PesinC1+}
  Let $f:M\circlearrowleft$ be a $C^{1+}$ diffeomorphism admitting an
  ergodic $f$-invariant probability measure which is hyperbolic and
  dominated. Then there exist constants
  $C,c,\theta, \delta_1>0, 0<\sigma<1$ and an integer $\ell\ge0$
  (depending only on $f$ and on the exponents of $\mu$) and measurable
  subsets $B^u,B^s$ with $\mu(B^*)>\theta, *=s,u$ such that
  \begin{enumerate}
  \item each $x\in B^s$ admits a stable manifold
    $\Delta=W^s_x(\delta_1)$ with inner radius at least $\delta_1$
    satisfying $\dist_{f^i\Delta}(f^iy,f^iz) \le
    \sigma^{i/2}\dist_\Delta(y,z)$ for all $y,z\in\Delta$ and all
    $i\in\ZZ^+$;
  \item each $x\in B^u$ admits an unstable manifold
    $\Delta=W^u_x(\delta_1)$ with inner radius at least $\delta_1$
    satisfying $\dist_{f^{-i}\Delta}(f^{-i}y,f^{-i}z) \le
    \sigma^{i/2}\dist_\Delta(y,z)$ for all $y,z\in\Delta$ and $i\in\ZZ^+$;
  \item $B:=B^u\cap f^{-\ell}B^s$ has positive $\mu$-measure and
    every $x\in B$ admits also a stable manifold $W^s_x(c)$ with
    inner radius $c$ and satisfying
    $\dist_{f^i\Delta}(f^iy,f^iz) \le C\sigma^{i/2}\dist_\Delta(y,z)$
    for all $y,z\in\Delta$ and $i\in\ZZ^+$.
  \end{enumerate}
  Moreover, the lamination $\cF^s:=\{W^s_x(\delta_1): x\in B^s\}$ is a
  continuous family of embedded disks which forms an
  absolutely continuous lamination, whose holonomy between $cu$-disks
  admits Jacobian bounded from above and from below away from zero.
\end{maintheorem}

For the meaning of absolute continuity and Jacobian of the holonomy 
along the stable leaves, see e.g.~\cite[Chapter 8]{BarPes2007} and
also Section~\ref{sec:long-unstable-leaves}.

\begin{remark}[comparison with Pesin's Non-Uniform Hyperbolic Theory]
  \label{rmk:Pesin}
  In the setting of the Non-Uniform Hyperbolic Theory of
  Pesin~\cite{BarPes2007} for $C^{1+}$ diffeomorphisms, or for
  hyperbolic and dominated probability measures for $C^1$
  diffeomorphisms, as considered by Abdenur et al. in~\cite[Theorem
  3.11, Section 8]{Abdenur2011}, we \emph{neither have a uniform
    contraction rate} on a neighborhood of uniform radius provided by
  the hyperbolic times; \emph{nor a global control of the curvature}
  of (un)stable disks.

  In particular, this means that the positive measure subset $B$,
  obtained from the \emph{coherent blocks} $B^u, B^s$ (see
  Section~\ref{sec:schedul-omega-limits} and~\cite{pinheiro20}), has
  stronger features than the \emph{hyperbolic blocks} from the
  Non-Uniform Hyperbolic Theory of Pesin\footnote{Even though coherent
    blocks cannot be enlarged to almost full measure.}. The discussion
  of \emph{effective hyperbolicity} by Climenhaga and Pesin
  in~\cite{ClmPe2016} is another example of the stronger features
  provided by hyperbolic times when coupled with non-zero Lyapunov
  exponents.
\end{remark}

% Similarly to the recent results from~\cite{LuoYang24} for $C^\infty$
% dynamics and not necessarily dominated measures, we obtain the
% following.

% \begin{maincorollary}[entropy and long (un)stable manifolds]
%   \label{mcor:hWus}
%   In the setting of Theorem~\ref{mthm:PesinC1+}, let $DH(f)$ be the
%   collection of all $f$-invariant ergodic hyperbolic and (non-atomic)
%   dominated probability measures; and $L^*(\delta)$ the subset of
%   points admitting an (un)stable manifold with inner radius larger
%   than $\delta>0$, where $*\in\{u,s\}$. Then, for every small enough
%   $\delta,\theta>0$, we can find $h>0$ so that if $\mu\in DH(f)$ is
%   such that $h_\mu(f)>h$, then
%   $\min\{ \mu(L^s(\delta)) , \mu(L^u(\delta)) \} > \theta$.
% \end{maincorollary}

\subsection{Organization of the text, comments and conjectures}
\label{sec:comments-conjectures}

We present examples of application of the main result to polynomial
mixing and robust exponential mixing for ergodic physical/SRB measures
for diffeomorphism with a dominated splitting, in the next
Section~\ref{sec:auxiliary-results}. In
Section~\ref{sec:hyperb-dominat-measu}, we present the main tools used
in the proofs, mainly from the recent book~\cite{Alves2020b} by Alves
and papers \cite{Pinho2011,pinheiro20} by one of the coauthors, and
references therein.

We provide a proof of Theorem~\ref{mthm:PesinC1+} in
Section~\ref{sec:long-unstable-leaves}.  In
Section~\ref{sec:existence-physic-mea-1} we present a proof
Theorem~\ref{mthm:NUHypGMY} together with most of
Corollary~\ref{mcor:abv}.  In Section~\ref{sec:speed-mixing-from}, we
deduce the statement of Theorem~\ref{mthm:Gibbsmix}.
Finally, in Section~\ref{sec:geometr-basin-coinci}
we deduce the statement of Theorem~\ref{mthm:basin} and the basin
claim of Corollary~\ref{mcor:abv}.

In the rest of this section we comment and conjecture possible
extensions of our results.

\subsubsection{Comments and conjectures}
\label{sec:comments-conjectures-1}

In all the previous main statements, we may replace the assumption on
the existence of \emph{dominated splitting} by the assumption of
existence of a \emph{$Df$-invariant and H\"older-continuous splitting}
$T_A M = E^{cs}_A\oplus E^{cu}_A$ and keep the same results --- it is
enough to follow the arguments from Cao, Mi and Yang~\cite{CMY2017}.

\begin{remark}[the assumption of existence of
  an ergodic hyperbolic $cu$-Gibbs state is not superfluous]
  \label{rmk:nonsuperfl} 
  Indeed, we consider a pair of diffeomorphisms
  $f,g:\cS^1\times\DD\circlearrowleft$, where $f$ is the uniformly
  hyperbolic Smale solenoid map, see
  e.g.~\cite[Sec. 7.7]{robinson1999}); and $g$ its ``intermittent''
  modification~\cite[Sec. 4.6]{Alves2020b} from~\cite{AlPi08}. In both
  cases we have attractors (i.e. transitive attracting sets)
  $\Lambda_f, \Lambda_g$ with partially hyperbolic splitting
  $E^s\oplus E^{cu}$ and ergodic (in fact, mixing) hyperbolic
  $cu$-Gibbs states $\mu_f,\mu_g$ which are the unique physical
  measures, but $f$ is uniformly hyperbolic, while $g$ admits a fixed
  point $p\in\Lambda_g$ so that $Dg_p\mid E^{cu}$ is an
  isometry. Hence, for $F:=f^\ell\times g$, $\nu=\mu_f\times\mu_g$ is
  an ergodic $F$-invariant measure which is the unique physical
  measure and a $cu$-Gibbs state, where $\ell>1$ is such that
  expansion/contraction rates of $f^\ell$ are stronger than the ones
  of $g$. Thus, $F$ is nonuniformly hyperbolic on a full volume
  measure subset of $(\cS^1\times\DD)^2$ and the attractor
  $\Lambda:=\Lambda_f\times\Lambda_g$ admits the dominated splitting
  $T_\Lambda M=\big(E^s_f\oplus E^s_g\big)\oplus \big(E^{cu}_g\oplus
  E^{cu}_f\big)$. However, the ergodic $cu$-Gibbs state
  $\nu=\mu_f\times\delta_p$ is non-hyperbolic, with a zero Lyapunov
  exponent along the direction $E^{cu}_g$. This shows that \emph{even
    with a full volume of non-uniformly hyperbolic points and unique
    physical/SRB measure there can be ergodic $cu$-Gibbs states which
    are not hyperbolic}.
\end{remark}

Recent results from Alves-Dias-Luzzatto-Pinheiro~\cite{ADLP} and
Bourguet-Yang~\cite{burguetYang24} allow us to obtain $cu$-Gibbs
states (which become physical measures) with \emph{weak non-uniform
  expansion}
\begin{align}\label{eq:wNUE0}
  \liminf\nolimits_{n\nearrow\infty} S_n\phi^{cu}(x)/n <0
\end{align}
on a positive volume subset of points in the trapping region. In their
partially hyperbolic setting, this \emph{a fortiori} implies
non-uniform expansion~\eqref{eq:NUE0} and so all our results can be
restated using this weak form of non-uniform expansion on a partially
hyperbolic setting.

In addition, it is natural to consider \emph{weak non-uniform
  contraction}
\begin{align}\label{eq:wNUC}
  \liminf\nolimits_{n\nearrow\infty} S_n\phi^{cs}(x)/n <0
\end{align}
on a positive volume subset of the trapping region.

We note that Tahzibi in~\cite{tah04} used the non-uniform
contraction~\eqref{eq:NUC} to obtain the existence of long stable
leaves Lebesgue almost everywhere, which then enables one to apply the
``Hopf argument'' to prove the existence of physical measures. It is
then natural to propose the following.
\begin{conjecture}
  \label{conj:wNUHyphys}
  Every attracting set with a dominated splitting with both weak
  non-uniform expansion~\eqref{eq:wNUE0} and weak non-uniform
  contraction~\eqref{eq:wNUC} admits a physical measure.
\end{conjecture}
If this holds true, then Theorem~\ref{mthm:Gibbsmix} applies to this
physical measure. 

We present in Subsection~\ref{sec:dominat-splitt-slowe} a non-robust
class of examples with polynomial rates of mixing for their physical
measures. It is natural to pose the following.

\begin{conjecture}
  \label{conj:robustpoly}
  There are examples of $C^r$ open subsets of diffeomorphisms
  ($r\ge1$) with dominated splitting together with non-uniform
  expansion and non-uniform contraction, without neither uniformly
  expanding nor contracting subbundles, having mixing physical
  measures which do not mix exponentially.
\end{conjecture}

The dependence of the rate of mixing exclusively from the tail
set of hyperbolic times along the unstable direction seems to follow
from the existence of a cylinder, in the ambient space, with a full
volume subset of long stable leaves with uniform contraction
rate. Therefore we pose the following.

\begin{conjecture}
  \label{conj:dependstable}
  There are examples of smooth diffeomorphisms, with hyperbolic
  physical measures, whose stable leaves admit no cylinder where their
  size is uniform, on a full volume subset, and whose mixing rates
  depend on the tail of hyperbolic times along the stable direction,
  that is, the analogous subset to~\eqref{eq:cutail} with $\phi^{cs}$ in
  the place of $\phi^{cu}$.
\end{conjecture}

Since the relation between geometric and ergodic basins, obtained in
Theorem~\ref{mthm:basin}, was a corollary of the existence of a GMY
structure, we pose the following.

\begin{conjecture}
  \label{conj:basin}
  There are smooth diffeomorphims with hyperbolic $cu$-Gibbs states
  whose ergodic basins are essentially different from their geometric
  basins.
\end{conjecture}

\begin{remark}
  \label{rmk:basin}
  This conjecture is false if we consider only physical measures, as
  the following example of a ``figure 8'' attracting set shows; see
  Figure~\ref{fig:8attractor}.

  \begin{figure}[htpb]
  \centering
  \includegraphics[width=9cm]{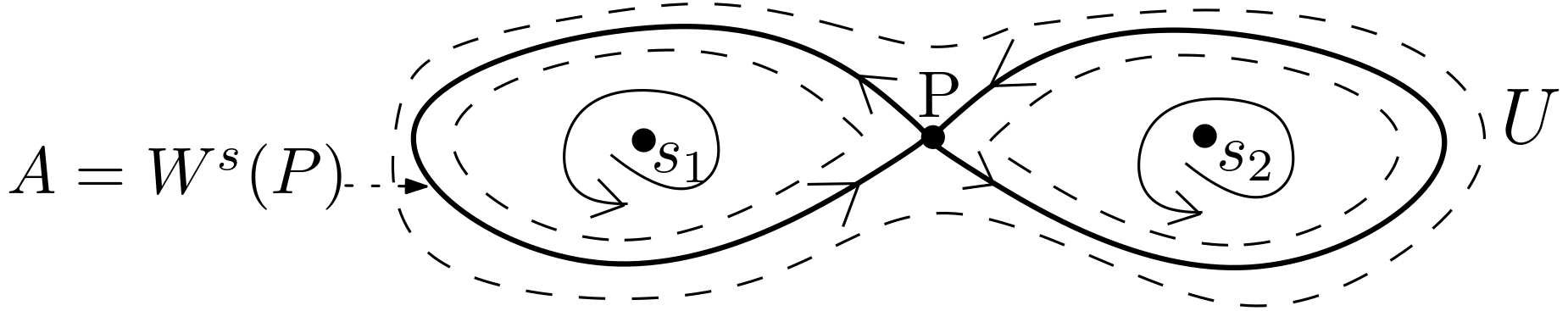}
  \caption{Sketch of the ``figure 8'' attracting set $A$ given by the
    double homoclinic connection $W^s(P)$ associated to tha hyperbolic
    saddle fixed point $P$ and an attracting neighborhood $U$.}
  \label{fig:8attractor}
\end{figure}

Indeed, we note that the only invariant measure supported on the
neighborhood of the invariant set $A$ is $\mu=\delta_P$ the Dirac mass
at the hyperbolic saddle fixed point $P$. Hence, this is also the only
accumulation point of the empirical measures
$\mu_n(x):=(1/n)S_n\vfi(x)$ for all $x$ in an open neighborhood $U$ of
$A$ as $n\nearrow\infty$. It follows that $B(\mu)\supset U$ and so
$\mu$ is an ergodic hyperbolic dominated and physical probability
measure.
  
However, the stable set $W^s(q)$ of each $q\in A$ coincides with
$W^s(P)=A$ and so $G(A)=A$ and $U\setminus G(A)=U\setminus A$ is an
open set, so the geometric and ergodic basin are essentially
different; see e.g.~\cite{CoYo2004}.
\end{remark}

\subsection*{Acknowledgements}

We thank the Mathematics and Statistics Institute of the Federal
University of Bahia (Brazil) for its support of basic research and
CNPq (Brazil) for partial financial support.

%%%%%%%%%%%%%%%%%%%%%%%%%%%%%%%%%%%%%%%%%%%%%%%%%%%%%%%%%%%%%%%%%
\section{Examples of application}
\label{sec:auxiliary-results}

Here present some examples of application of the main theorems. In
Subsection~\ref{sec:partially-hyperb-exa}, we consider partially
hyperbolic examples with uniformly expanding subbundle as particular
applications of the main results, obtaining exponential mixing for
physical measures. A non-robust class of examples with slower
(polynomial) rates of mixing is presented in
Subsection~\ref{sec:dominat-splitt-slowe}. A robust class of
exponential mixing for physical measures of partially hyperbolic and
non-uniformly hyperbolic diffeomorphisms without uniform invariant
subbundle is described in Subsection~\ref{sec:robust-non-uniform}.

\subsection{Partially hyperbolic examples}
\label{sec:partially-hyperb-exa}

We start with partially hyperbolic examples with spliting
$E^{cs}\oplus E^u$, where $E^u$ is uniformly expanding, enabling us to
more easily obtain ergodic physical/SRB measure which are also
$cu$-Gibbs states with exponential mixing, independently of the fine
asymptotic behavior along the center-stable direction.

The examples presented in the works of Bonatti-Viana~\cite{BoV00} and
Castro~\cite{Castro04,Castro02} provide robust families of $C^2$
diffeomorphisms with partially hyperbolic splitting admitting
physical/SRB ergodic probability measures. Since in this cases we have
uniform expansion along the unstable direction, we have non-uniform
expansion and the average expansion time function $h$ is constant on
the ergodic basin of the physical measures. We immediately obtain from
Theorem~\ref{mthm:NUHypGMY} and Corollary~\ref{mcor:GMYexp} the
following.

\begin{corollary}
  \label{cor:phypexpmix}
  Let $f:M\circlearrowleft$ be a $C^{1+\eta}$ diffeomorphism, for some
  $\eta\in(0,1]$, with a partially hyperbolic splitting
  $T_AM=E^{cs}_A\oplus E^{u}_A$ over an attracting set $A$ on a
  trapping region $U\subset M$, and an ergodic physical/SRB measure
  $\mu$. Then there exists a power $g=f^q$ for some $q\ge1$ such that
  there are $1\le p\le q$ invariant exponentially mixing probability
  measures $\nu_1,\dots,\nu_p$ so that $f_*\nu_i=\nu_{i+1}$ for
  $i=1,\dots, p-1$; $f_*\nu_p=\nu_1$, and
  $p\cdot\mu=\sum_{i=1}^p \nu_i$. More precisely, for each
  $1\le i\le p$ we can find $c>0$ so that $\eta$-H\"older observables
  $\vfi,\psi:M\to\RR$ admit $C>0$ for which
  $\corr_{\nu_i}(\vfi,\psi\circ f^{qn})\le Ce^{ - c n}$ for all
  $n\ge1$.
\end{corollary}

\begin{remark}\label{rmk:stablecontraction}
  We note that we have no condition on the ``average contraction
  function'' along the central-stable direction.
\end{remark}

\subsection{Dominated splitting and slower rates of mixing}
\label{sec:dominat-splitt-slowe}

We describe an example of a non-uniformly hyperbolic attractor with
dominated splitting with a unique physical measure which is
polynomially mixing independently of the eventual rates of convergence
along the center-stable direction.

We recall the construction of the solenoid with intermittency
from~\cite[Sec. 2.4]{AlPi08}.  Let $f: \sS^1\circlearrowleft$ be a map
of degree $d\ge 2$ with the following properties:
\begin{enumerate}
  \item[(i)] $f$ is $C^2$ on $\sS^1\setminus\{0\}$;
  \item[(ii)]  $f$ is $C^1$ on $\sS^1$ and $f'>1$ on $\sS^1\setminus\{0\}$;
  \item[(iii)] $f(0)=0$, $f'(0)=1$, and there is $\gamma>0$ such that
    $-xf''(x)\approx |x|^\gamma\quad\text{for all $x\neq 0$.}$
\end{enumerate}
Consider the solid torus $M=\sS^1\times\DD^2$, where $\DD^2$ is the
unit disk in $\CC$, and define $F: M\circlearrowleft$ by
$F(x,z):=\big(f(x), g(\theta,z)\big)$ where $g(\theta,z):=(z/10+e^{ix}/2).$

From~\cite[Sec. 5.1]{AlPi08} (cf.~\cite{You99}) we have that $f$
admits an absolutely continuous ergodic invariant probability measure
$\nu$ if, and only if, $\gamma<1$; and, moreover, $f$ is non-uniformly
expanding whose average expansion function $h$  satisfies
$\lambda(\{h>n\})\le C n^{-1/\gamma}$.  Since $F$ is conformal along
$\DD$ and uniformly contracting, we are in the setting of non-uniform
hyperbolicity and recover the results from~\cite{AlPi08}.

However, we can modify $g$ on a neighborhood of a periodic orbit to
obtain non-uniform contraction keeping the non-uniform expanding
structure of $F$. Indeed, since $f$ has degree two, then there exists
a period-two orbit $\{\theta_0,\theta_1:=f(\theta_0)\}$ for $f$ and
$F^2(\theta_0,z) = (\theta_0, g(f(\theta_0),g(\theta_0,z))) =
(\theta_0,g_2(\theta_0,z))$ where $z\mapsto g_2(\theta_0,z)$ is a
conformal $1/100$-contraction on $\DD$. Hence, there is a fixed point
$z_0\in\DD$ for this action so that
$F(\theta_0,z_0)=(\theta_1,g(\theta_0,z_0))=(\theta_1,z_1)$ and
$F(\theta_1, z_1)=(\theta_0,z_0)$.

We perform a $C^\infty$ modification of $g$ on small neighborhoods $V_0$ of 
$(\theta_0,z_0)$ and $V_1$ of $(\theta_1,z_1)$ so that the new
function $\wt{g}:M\to\sS^1$ keeps a conformal derivative and also, writing
$\wt{g_2}(\theta,z):= \wt{g}\big(f(\theta), \wt{g}(\theta,z)\big)$: 
\begin{enumerate}[(a)]
\item $D_2\wt{g_2}(\theta_0,z_0)=1$ and;
\item $D_2\wt{g_2}(\theta,z)<1$ for all $(\theta,z)\notin\{
  (\theta_0,z_0),(\theta_1,z_1)\}$.
\end{enumerate}
It is easy to see that $\phi^{cs}=\log\|D_2\wt{g_2}\|$, where
$E^{cs}_{(\theta,z)}\approx\RR^2$ is the tangent space $T_z\DD$ and
$\xi(\theta):=\max_{z\in\DD}\phi^{cs}(\theta,z)$, satisfies 
$\int\xi(\theta) \,d\nu(\theta) <0$. Since $\nu\times\leb_{\DD}$, with
$\leb_{\DD}$ the Lebesgue measure on the disk $\DD$, is equivalent to
Lebesgue measure $\leb$ on $M$, then non-uniform
contraction~\eqref{eq:NUC} for
$\wt{F}(\theta,z):=\big(f(\theta),\wt{g}(\theta,z)\big)$
follows. Indeed, % writing $(\theta_i,z_i)=F^i(\theta,z), i\ge1$
for $\nu$-a.e. $\theta\in\sS^1$ and each $z\in\DD$, we get a point
$x=(\theta,z)\in M$ satisfying
\begin{align*}
  \limsup\nolimits_{n\nearrow\infty}S_n^F\phi^{cs}(x)/n
  \le
  \limsup\nolimits_{n\nearrow\infty}S_n^f\xi(\theta)
  =
  \int\xi\,d\nu<0.
\end{align*}
Moreover $\phi^{cu}$, with $E^{cu}=T\sS^1$, is non-uniform expanding
since $D\wt{F}\mid E^{cu}=D_1\wt{F}=Df\circ\pi$ where $\pi:M\to\sS^1$
is the canonical projection into the first coordinate.

In addition, the map $\wt{F}$ is $C^\infty$ and the
$D\wt{F}$-invariant spitting $TM=E^{cs}\oplus E^{cu}$ is dominated
because (recall that $D_2g$ is conformal)
\begin{align*}
  \frac{\|D\wt{F}\mid E^{cs}_{(\theta,z)}\|}{\|D\wt{F}\mid
  E^{cu}_{(\theta,z)}\|}
  =
  \frac{D_2\wt{g}(\theta,z)}{Df(\theta)}
  \le
  \begin{cases}
    Df(\theta)^{-1}, & \theta\neq0
    \\
    D_2\wt{g}(0,z), & \theta=0
  \end{cases}
\end{align*}
is a continuous function $M\to\RR$ strictly smaller than $1$.

Since $\wt{F}$ is transitive on $A$ as a consequence of the
transitivity of $f$, therefore the attractor
$A=\cap_{n\ge0}\wt{F}^n(M)$ admits a unique ergodic physical/SRB
measure which is also a $cu$-Gibbs state $\mu$.

We can now follow the construction presented in
Section~\ref{sec:existence-physic-mea-1} to check that we are in the
case (1) of the statement of Theorem~\ref{mthm:Gibbsmix}, with $q=1$,
obtaining polynomial mixing for this attractor. Indeed, since $f$ is
topologically mixing, then $\wt{F}$ is topologically mixing on $A$ and
then we can take the power $q=1$ to obtain mixing for the measure
$\mu$ with respect to the action of $\wt{F}$.

\begin{remark}\label{rmk:nonrobustpoly}
  This example is not robust since the non-uniform expansion depends on
  the tangency of the graph of the function $f$ to the diagonal; see
  e.g.~\cite{ArTah,ArTah2}.
\end{remark}

\subsection{Robust  example of exponential mixing for physical
  measures without uniform invariant subbundle}
\label{sec:robust-non-uniform}
The $C^1$ open classes of transitive non-Anosov diffeomorphisms
presented in \cite[Section~6]{BoV00}, as well as other robust examples
from \cite{Man78}, and also in~\cite{ABV00} and~\cite{tah04} are
constructed in a similar way.

\subsubsection{General description of the geometric properties}
\label{sec:general-descript-geo}

We assume that we start with some Anosov diffeomorphism $\hat f$ on
the $d$-dimensional torus $M=\TT^d$, $d\ge 3$ with a decomposition of
the tangent fiber bundle $TM=E^{uu}\oplus E^{ss}$. Let $W$ be an open
subset in $M$ and let us assume that that $f$ is a $C^1$ close
diffeomorphism satisfying
 \begin{enumerate}[(A)]
 \item the tangent bundle decomposes $TM=E^{cs}\oplus E^{cu}$ into a
   dominated splitting and $f$ admits invariant cone fields $C^{cu}$
   and $C^{cs}$, with small width $a>0$ and containing, respectively,
   $E^{cu}$ and $E^{cs}$;
 \item $f$ is \emph{volume hyperbolic}: there is $\sigma_1>1$ so
   that
   $$ |\det(Df\vert T_x\D^{cu})| > \sigma_1 \quad\mbox{and}\quad
   |\det(Df\vert T_x\D^{cs})| < \sigma_1^{-1} $$ for any $x\in M$ and
   any disks $\D^{cu}$, $\D^{cs}$ tangent to $C^{cu}$, $C^{cs}$,
   respectively.
 \item $f$ is $C^1$-close to $\hat f$ in the complement of $W$,
   so that there exists $\sigma_2<1$ satisfying
 $$
 \|(Df \vert T_x \D^{cu})^{-1}\| < \sigma_2
\quad\mbox{and}\quad \|Df \vert T_x \D^{cs}\| < \sigma_2
 $$
 for any $x\in (M\setminus W)$ and any disks $\D^{cu}$, $\D^{cs}$
 tangent to $C^{cu}$, $C^{cs}$, respectively.
 \item there exist some small $\delta_0>0$ satisfying
 $$
 \|(Df \vert T_x\D^{cu})^{-1}\| < 1+\delta_0\quad\mbox{and}\quad \|Df
 \vert T_x \D^{cs}\| < 1+\delta_0
 $$
 for any $x\in W$ and any disks $\D^{cu}$ and $\D^{cs}$ tangent to
 $C^{cu}$ and $C^{cs}$, respectively.
\end{enumerate}

\subsubsection{Robust non-uniformy hyperbolic example}
\label{sec:non-partially-hyperb}

From \cite[Theorem C]{BoV00}, \cite[Appendix]{ABV00} together with
Tahzibi~\cite{tah04}, performing a small perturbation along the
central-stable and center-unstable direction of the initial Anosov
diffeomorphism $\hat f:\TT^d\circlearrowleft$ with $d\ge4$ on the
region $W$, provides the following; see also~\cite[Section
7.1.4]{BDV2004}.

\begin{figure}[htpb]
  \centering
  \includegraphics[width=9cm]{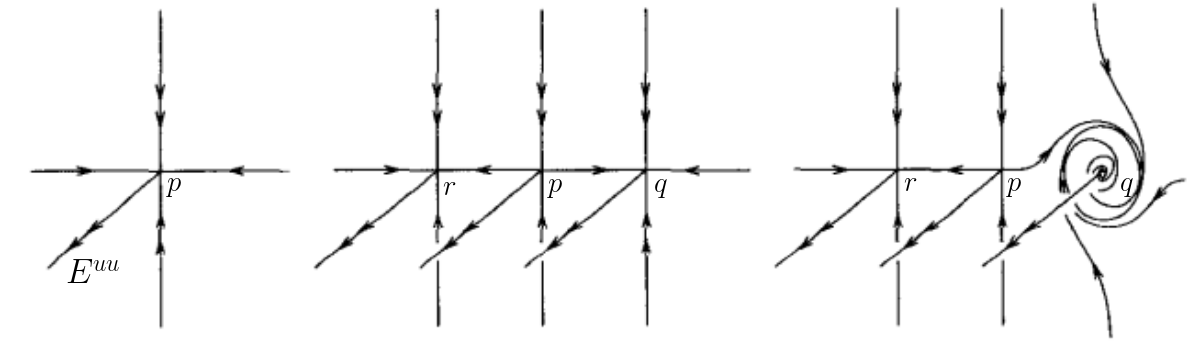}
  \caption{Sketch of the deformation of the linear Anosov
    diffeomorphism around the hyperbolic fixed point $p$ with stable
    index $s$ in the left hand side. In the center figure two new
    saddles appear with the same stable index $s$ while the
    stable index of $p$ becomes $s-1$. In the right hand side,
    the saddle $q$ becomes an attracting center along the stable
    direction. The strong unstable direction $E^{uu}$ depicted above
    has dimension $u\ge2$.}
  \label{fig:deformAnosov}
\end{figure}

\begin{theorem}
  \label{thm:NUEnotPH}
  There exists a $C^2$ neighborhood $\cV$ of $f$ and $c_u,c_s>0$ such
  that all diffeomorphisms $g\in\cV$ are topologically mixing with a
  non-uniformly hyperbolic dominated splitting
  $T\TT^d=E^{cs}\oplus E^{cu}$. Moreover, $g$ admits no other
  invariant subbundle, and $\cV$ contains an open subset of the space
  of $C^2$ volume preserving diffeomorphisms of $\TT^d$.

  In addition, there exists a periodic point $p\in M\setminus W$ whose
  stable $W^s_p$ and unstable $W^u_p$ manifolds are dense for each
  $g\in\cV$; and there exists a unique physical/SRB measure $\mu_g$,
  which is also the unique $cu$-Gibbs state with full basin
  $\leb(M\setminus B(\mu_g))=0$. %In addition, $\mu_g(W)>0$.
\end{theorem}

\begin{proof}
  This is the main result of Tahzibi in~\cite{tah04}, which proves all
  statements. % except the last on positive mass of $W$.
  The deformation of the Anosov diffeomorphism $\hat f$ on $\TT^4$
  starting with a hyperbolic decomposition $TM=E^{ss}\oplus E^{cu}$
  with $s=\dim E^{ss}=\dim E^{uu}=u=2$, can be described as
  follows\footnote{This can easily be extended to any dimension
    $d=s+u\ge4$ with $s,u\ge2$; see~\cite{tah04}
    for more details.}; see Figure~\ref{fig:deformAnosov}.

  We fix a small neighborhood $W$ of a fixed point $p$ of $\hat f$ (or
  of a power $\hat f^k$ if needed) and take a one-parameter
  family\footnote{For more details on the construction of this family,
    see~\cite[Section 6.4]{BoV00}.}  of diffeomorphisms
  $(f_t)_{t\in[0,1]}$ so that, as first stage:
    \begin{enumerate}[(I)]
    \item the point $p$ is fixed for every $f_t$;
    \item the weakest contracting eigenvalue of $Df_t(p)$ increases as
      $t$ increases from $0$;
    \item at some $0<t=t_0<1$ this eigenvalue becomes equal to $1$,
      and the stable index (dimension of the stable bundle) of $p$
      changes from $2$ to $1$;
    \item in the process, for $t=t_1\in(t_0,1)$, new fixed saddle
      points $r,q$, with stable index $2$, are created in the
      neighbourhood of $p$.
    \end{enumerate}
    At this stage, for $t_1$ close to $t_0$, if we set $g_0=f_{t_1}$,
    then $g_0$ admits a partially hyperbolic $Dg_0$-invariant
    splitting $TM=E^s\oplus E^{cu}$ so that $E^{cu}$ is close to
    $E^{uu}$ and $E^s$ close to the original stable bundle of
    $\hat f$; and also $E^{cu}$ is non-uniformly expanding for a
    certain rate $c_u>0$. For details, see e.g.~\cite[Appendix]{ABV00}
    or~\cite{tah04}.
    
    We consider a small neighborhood $V_q$ of
    the saddle $q$ such that $V_q\subset W$ but does not contain
    $p,r$. Then proceed to the second stage, modifying $g_0$ in this
    neighborhood obtaining a one-parameter family $g_s$ of
    diffeomorphisms so that
    \begin{enumerate}[(i)]
    \item $q$ is a fixed point of every $g_s$;
    \item the contracting eigenvalues of $Dg_s(q)$ become equal, and
      then complex conjugate, as $s$ becomes larger than some small
      $s_0>0$.
    \end{enumerate}
    Let $h=g_{s_1}$ for some $s_1>s_0$ close to $s_0$. We can perform
    these changes keeping the stable foliation of $\hat f$ still
    $h$-invariant so that any sufficiently thin cone field around the
    stable foliation of $\hat f$ is a centre-stable cone field for
    $h$; and also ensure that there exists a sufficiently thin
    center-unstable cone field around the initial unstable direction.

    To complete the construction, we repeat the deformation steps
    (i)-(ii) outlined above starting from the diffeomorphism $h$ for a
    small neighborhood $V_r$ around the saddle $r$, in the place of
    the saddle $q$, where $V_r$ does not contain $p,q$ but is
    contained in $W$; and \emph{the expanding eigenvalues are used in
      the place of the contracting eigenvalues} in step (ii). This
    diffeomorphism $f$ admits a $Df$-invariant
    dominated decomposition $TM=E^{cs}\oplus E^{cu}$, with $E^{cs}$
    non-uniformly contracting for some rate $c_s>0$ and $E^{cu}$ still
    non-uniformly expanding.

    This provides us a with the diffeomorphim $f$ and the $C^2$
    neighborhood $\cV$ in the statement of Theorem~\ref{thm:NUEnotPH},
    as shown in~\cite{tah04}.
    % To conclude, we need to show that the physical/SRB measure $\mu_g$
    % gives positive mass to $W$. We observe that the denseness of the
    % stable and unstable manifolds of any given fixed periodic point
    % $w\in M\setminus W$ ensures, by the $\lambda$-Lemma, that any
    % center-unstable disk in the support of $\mu_g$ will eventually
    % intersect the open subset $W$. Then $\mu_g(W)>0$.
\end{proof}

\subsubsection{Robust exponential mixing}
\label{sec:robust-exponent-mixi}

The reader should recall the expansion time function $h$ from
Subsection~\ref{sec:exponent-mixing-phys}.

\begin{proposition}\label{pr:robustexpmix}
  Every $f\in\cV$ is such that every $cu$-disk $\gamma\subset M$
  admits a subset $H\subset\gamma$ with a full $\leb_\gamma$-measure
  where $f$ is non-uniformly hyperbolic and $\leb_\gamma(h>n)$ decays
  exponentially fast to $0$ with $n$.
\end{proposition}

\begin{proof}
  This follows from the arguments in~\cite[Appendix]{ABV00} or, with a
  more detailed presentation, from~\cite[Proposition
  7.32]{Alves2020b}.
\end{proof}

Proposition~\ref{pr:robustexpmix} together with
Corollary~\ref{mcor:GMYexp} ensures that the family $\cV$ is a $C^2$
\emph{robust family of non-uniformly hyperbolic exponentially mixing
  diffeomorphisms without any uniformly contracting or expanding
  invariant subbundle}.

\section{Auxiliary results}
\label{sec:hyperb-dominat-measu}

Here we present the tools used in the proofs of the main results.
From now on we assume that $f$ is a $C^{1+}$ diffeomorphism with a
compact invariant attracting subset $A$ with trapping region $U$
admitting a dominated splitting which is non-uniformly hyperbolic for
a $\leb$-positive subset $H$ of $U$.  We assume without loss that the
splitting has been continuously extended to the open neighborhood $U$
of $A$ and that all constructions are performed in this
(relatively compact) neighborhood.

% \subsection{Birkhoff generic points}
% \label{sec:birkhoff-generic-poi}

% The Borel $\sigma$-algebra $\B$ of $M$ is countably generated by the subset
% $\G$ of balls of rational radius centered on points of a dense subset
% of $M$. This enables us to obtain:

% \begin{lemma}[Full measure of \emph{Birkhoff generic points}]
%   \label{le:Birkhoffgeneric}
%   Given an $f$-ergodic probability measure $\mu$ on the complete
%   separable metric space $M$, then there exists a full $\mu$-measure
%   subset $Y=Y(\mu)$ of \emph{Birkhoff generic points}, that is, those
%   points $x$ satisfying\footnote{As usual $1_E$ is the indicator function
%     of the set $E$: $1_E(x)=1$ if $x\in E$ and $1_E(x)=0$ otherwise.}
%   \begin{align}\label{eq:Birkhoffgeneric}
%     \tau(E,x):=\lim\nolimits_{n\nearrow\infty}S_n 1_E(x)/n=\mu(E)
%     \text{  for all  } 
%     E\in\B.
%   \end{align}
% \end{lemma}

% \begin{proof}
%   Apply the Ergodic Decomposition Theorem~\cite[Section 5]{ViOl16} to
%   $\mu$ and observe that the unique component must satisfy the above
%   statement for the mean sojourn time $\tau(E,x)$ with respect to
%   $\mu$-a.e. $x\in M$ of any element $E$ of $\B$: cf. the proof of
%   Ergodic Decomposition on complete separable metric spaces
%   at~\cite[Section 5.1.4]{ViOl16}. We let $Y=Y(\mu)$ be this full
%   measure subset.
% \end{proof}

\subsection{Consequences of the existence of dominated splitting}
\label{sec:conseq-existence-dom}

From the existence of dominated splitting, it is a standard
fact\footnote{See e.g.~\cite[Theorem 5.5]{HPS77} or the statement of
  \cite[Lemma 4.4]{ArMel17}.}  that there are continuous families
$(W^{*}_x)_{x\in M}$ of $C^1$ embedded $*$-disks such that
$T_xW^*_x=E^*_x$ for $*\in\{cs,cu\}$ and locally invariant, i.e. for
each $0<\epsilon<\epsilon_0$ and all $x\in A$ there exists $\delta>0$
such that
\begin{align}\label{eq:localinv}
  f^{-1}(W^{cs}_x(\epsilon_0))\cap B_\delta(f^{-1}x)\subset
  W^{cs}_x(\epsilon) \qand  f(W^{cu}_x(\epsilon_0))\cap B_\delta(fx)\subset
  W^{cu}_x(\epsilon),
\end{align}
where $W^*_x(\epsilon)$ is the $\epsilon$-ball in $W^*_x$ around $x$.

Given a $cu$-disk $\Sigma$, then $f(\Sigma)$ is also tangent to the
centre-unstable cone field by the domination property. The tangent
bundle of $\Sigma$ is said to be {\em H\"older continuous} if
$x \mapsto T_x\Sigma$ is a H\"older continuous section from $\Sigma$
to the Grassman bundle of $M$. In other words, at every $x\in\Sigma$
we can find a neighborhood $V$ where the $V\cap\Sigma$ is a graph of a
H\"older-$C^1$ function $\psi_x:E^{cu}_x\to E^{cs}_x$.  We define
\begin{align}
  \label{eq:kappa}
  \kappa(\Sigma)
  :=
  \inf\{C>0:\text{the tangent bundle of $\Sigma$
  is $(C,\zeta)$-H\"older}\},
\end{align}
where $\zeta>0$ is so that
$\|Df^n\mid E^{cs}_x\|\cdot\|(Df^n\mid
E^{cu}_{f^nx})^{-1}\|^{1+\zeta}$ still tends to zero when
$n\nearrow\infty$ for $x\in A$. The next result contains the
information needed on the H\"older control of the tangent direction.

%\margem{reobtain w/cont. splitting}

\begin{proposition}{\cite[Corollary 2.4]{ABV00}}
  \label{pr:curvature} There exists $C_1>0$ such that, given any $C^1$
  $cu$-disk $\Sigma\subset U$ such that
  $\Sigma\cap A\neq\emptyset$, then there exists $n_0\ge 1$ such
  that $\kappa(f^n(\Sigma)) \le C_1$ for every $n\ge n_0$. Moreover
  \begin{enumerate}
  \item if $\kappa(\Sigma) \le C_1$, then $\kappa(f^n(\Sigma)) \le C_1$ for
    every $n\ge 1$;
  \item if $\Sigma$ and $n$ are as above, then the functions
$$
J_k: f^k(\Sigma)\ni x \mapsto \log |\det \big(Df \mid T_x f^k(\Sigma)\big)|,
\quad\text{$0\le k \le n$},
$$
are $(L_1,\zeta)$-H\"older continuous with $L_1>0$ depending only on
$C_1$ and $f$.
\end{enumerate}
\end{proposition}

\subsection{Hyperbolic times and center-unstable pre-disks}
\label{sec:hyperbolic-times-pre}

We derive uniform expansion and bounded distortion estimates from the
non-uniform expansion assumption in the centre-unstable direction.

\begin{figure}[htpb]
    \includegraphics[width=9cm]{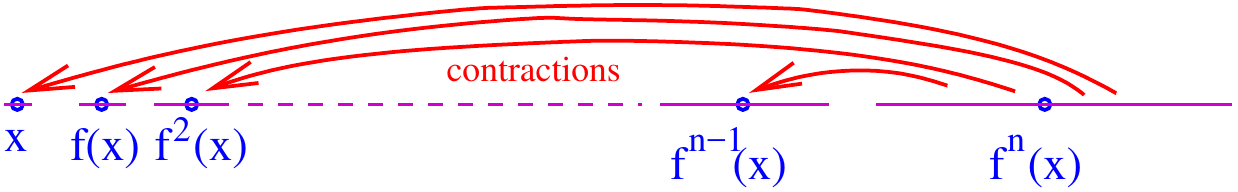}
    \caption{\label{fig:hyptime}Backward contractions at hyperbolic
      times.}
  \end{figure}

We say that $n$ is a $\sigma$-\emph{hyperbolic time} for $x\in U$ if
$0<\sigma<1$ and
\begin{align*}
  S_k\phi^{cu}(f^{n-k+1}x) \le k \log\sigma, \quad
  0\le k<n.
\end{align*}
In this case, $Df^{-k} \mid E^{cu}_{f^{n}(x)}$ is a contraction for
every $1\le k \le n$; see Figure~\ref{fig:hyptime}.
  
If $a>0$ is sufficiently small
% (recall the definition of the
%cone fields)
and we choose $0<\delta_1<\epsilon_0/2$ then, by
continuity
\begin{align}
  \label{eq:delta1} 
  \|Df(y)u\|\le \sigma^{-1/4}\|Df\mid E^{cs}_x\|\,\|u\|
  \;\,\&\,\;
  \|Df^{-1}(f(y)) v \| \le \sigma^{-1/4}
  \|(Df|E^{cu}_x)^{-1}\|\,\|v\|,
\end{align}
whenever $x,y\in M$, $d(x,y)\le \delta_1$, $u\in C^{cs}_a(y)$ and
$v\in C^{cu}_a(y)$.

Given any disk $\Delta\subset M$, we use $\dist_\Delta(x,y)$ to denote
the distance between $x,y\in \Delta$, measured along $\Delta$. The
distance from a point $x\in \Delta$ to the boundary of $\Delta$ is
$\dist_\Delta(x,\partial \Delta)= \inf_{y\in\partial
  \Delta}\dist_\Delta(x,y)$. The following has been proved in
\cite[Lemma 2.7]{ABV00}; see \cite[Lemma 4.2]{AlPi08a} for a detailed
proof.

\begin{lemma}[Pre-disks at hyperbolic times]
  \label{le:contraction} Let $0<\delta<\delta_1<\epsilon_0$,
  $0<\sigma<1$ and $\Delta\subset U$ be a $cu$-disk of radius
  $\delta$.
  % Let $\Delta$ be a $C^1$ centre-unstable disk contained in $U$.
  Then, there is $n_0\ge1$ such that for $x\in \Delta$ with
  $\dist_\Delta(x,\partial \Delta)\ge \delta/2$ and $n \ge n_0$ a
  $\sigma$-hyperbolic time for $x$ there is a neighborhood $W_n=W_n(x)$ of
  $x$ in $\Delta$ such that:
  \begin{enumerate}
  \item $f^{n}$ maps $W_n$ diffeomorphically onto a $cu$-disk of
    radius $\delta_1$ around $f^{n}(x)$;
  \item for every $1\le k \le n$ and $y, z\in W_n$:
    $$ \dist_{f^{n-k}(W_n)}(f^{n-k}(y),f^{n-k}(z)) \le
    (\sigma^{1/2})^k\dist_{f^n(W_n)}(f^{n}(y),f^{n}(z)).$$
  \end{enumerate}
\end{lemma}

\begin{remark}[Pre-disks and dynamical balls]
  \label{rmk:prediskdynball}
  Hence, each $y\in W_n$ has $n$ as a $\sigma^{1/2}$-hyperbolic time
  and $W_n$ is the $(n+1,\delta_1)$-\emph{dynamical ball around $x$ in
    $\Delta$}. That is, we have $W_n=\Delta\cap B(x,n+1,\delta_1)$,
  where we write, as usual,
  $B(x,n,\delta_1):=\{z\in M: d(f^iz,f^ix)<\delta_1, i=0,\dots,n-1\}$
  for the $(n,\delta_1)$-dynamical ball around $x$ in $M$.

  Moreover, from~\eqref{eq:delta1}, we have that any $cu$-disk
  $\gamma$ on $B(x,n+1,\delta_1)$ has $n$ as a
  $\sigma^{1/2}$-hyperbolic time for each $z\in\gamma$.
\end{remark}

We call the sets \( W_n \) \emph{hyperbolic pre-disks} and their
images \( f^{n}(W_n) \) \emph{hyperbolic disks}, which are indeed
centre-unstable balls of radius \( \delta_1 \).  The following is a
consequence of Proposition~\ref{pr:curvature} and
Lemma~\ref{le:contraction} above exactly as in the proof of
\cite[Proposition 2.8]{ABV00}.

\begin{corollary}[Bounded distortion]
  \label{cor:distortion} There exists $C_2>1$ such that given a disk
  $\Delta$ as in Lemma~\ref{le:contraction} with
  $\kappa(\Delta) \le C_1$, and given any hyperbolic pre-ball
  $W_n\subset \Delta$ with $n\ge n_0$, then 
$$
\log \frac{|\det Df^{n} \mid T_y \Delta|} {|\det Df^{n} \mid T_z
  \Delta|} \le C_2 \dist_{f^n(W_n)}(f^{n}(y),f^{n}(y))^\zeta,
\text{  for all  } y,z\in W_n.
$$
% Moreover, we may take $C_2=\exp(L_1\delta_1/(1-\sigma^{\zeta/2}))$.
\end{corollary}

The next result states the existence of hyperbolic times with positive
asymptotic frequency for points satisfying~\eqref{eq:NUE0} and its proof
can be found in \cite[Lemma 3.1, Corollary 3.2]{ABV00}.

\begin{proposition}[Positive frequency of hyperbolic times]
  \label{pr:hyperbolic2}
  For every $x \in U$ with $S_n\phi^{cu}(x) \le -c_un$ there exist
  $\sigma_u$-hyperbolic times $1 \le n_1 < \cdots < n_l \le n$ for
  \( x \) with $l\ge\theta_u n$ and $\sigma_u:=e^{-7c_u/8}$, where
  \( \theta_u:=c_u/(8\bar\phi^{cu}-7c_u) \) and
  $\bar\phi^{cu}:=\sup\{-\phi^{cu}(x): x\in U\}$.
\end{proposition}

% The following result ensures that forward invariant subsets which
% intersect the nonuniformly expanding set with positive volume must
% essentially contain centre-unstable disks.

% \begin{proposition}{\cite[Proposition 7.6]{Alves2020b}}
%   For each measurable subset $A\subset U$ such that $f(A)\subset A$
%   and $\leb(A\cap H)>0$ there are a $cu$-disk $D\subset U$ with
%   $\kappa(D) \le C_1$, integers $1 \le n_1 < n_2 < \ldots$ and
%   hyperbolic pre-balls $W_{n_k} (x_k )\subset D$ for each $k\ge1$, so
%   that for $\Sigma_k = f^{n_k} (W_{n_k} (x_k))$ we have
%   $\lim \frac{\leb_{\Sigma_k}(A\cap\Sigma_k)}{\leb_{\Sigma_k}(\Sigma_k)}=1$.
% \end{proposition}

\subsection{Reverse/Inverse hyperbolic times and center-stable
  pre-disks}
\label{sec:reverse-hyperb-times}

By assumption~\eqref{eq:NUC}, we have $c_s>0$ and a strictly
increasing sequence $m_i\nearrow\infty$ so that
$ S_{m_i}\phi^{cs}(x)<-c_sm_i$ as $i\nearrow\infty$.

\begin{figure}[htpb]
  \centering \includegraphics[width=9cm]{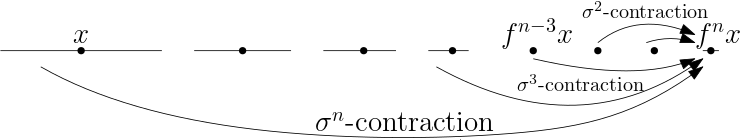}
  \\
  \includegraphics[width=9cm]{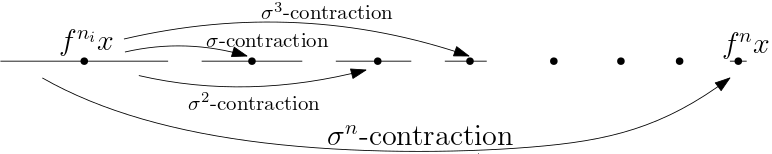}
  \caption{\label{fig:revhyptime}Forward contractions at inverse
    hyperbolic times above versus forward contractions at reverse
    hyperbolic times below.}
\end{figure}

Analogously to hyperbolic times  in the center-unstable direction, we
say that $n\ge1$ is a $\sigma$-\emph{inverse hyperbolic time} if
$0<\sigma<1$ and
\begin{align*}
 S_k\phi^{cs}(f^{n-k}x) \le k\log \sigma, \quad 0< k \le n;
\end{align*}
and that $n\ge0$ is a $\sigma$-\emph{reverse hyperbolic time} with
respect to $m>n$ if
\begin{align*}
  S_k \phi^{cs}(f^nx)\le k\log\sigma, \quad 0<k\le m-n.
\end{align*}
In Figure~\ref{fig:revhyptime} we depict the difference between
inverse and reverse hyperbolic times.

To ensure the existence of these times in our setting we use the
following.

\begin{lemma}[Pliss Lemma; see e.g. Chapter IV.11 in~\cite{Man87}]
  \label{le:pliss}
  Let $L\ge c_2 > c_1 >0$ and $\theta={(c_2-c_1)}/{(L-c_1)}$. Given
  real numbers $a_1,\ldots,a_N$ satisfying
  $ \sum_{j=1}^N a_j \ge c_2 N$ and $a_j\le L$ for $1\le j\le N$,
  there are $\ell>\theta N$ and $1<n_1<\ldots<n_\ell\le N$ such that
  $ \sum_{j=n+1}^{n_i} a_j \ge c_1\cdot(n_i-n) \;\;\mbox{for each}\;\;
  0\le n < n_i, \; i=1,\ldots,\ell.  $
\end{lemma}
We set $c_2=c_s$, $c_1=7c_2/8$,
$L=\bar\phi^{cs}:=\sup\{x\in U : -\phi^{cs}(x)\}$ and  we define 
\begin{enumerate}[(a)]
\item either $a_j=-\log\|Df\mid E^{cs}_{f^jx}\|$ ;
\item or $a_j=-\log\|Df\mid E^{cs}_{f^{m_i-j}x}\|$;
\end{enumerate}
for $1< j\le m_i$.  \emph{We note that we are inverting the summation
  order in the second case.}

Then, for $\theta_s=c_s/(8\bar\phi^{cs}-7c_s)>0$ and $N=m_i$, Pliss
Lemma~\ref{le:pliss} ensures that there are $\ell>\theta_s N$ and
$1<n_1<\dots<n_\ell\le m_i$ such that for each $k=1,\ldots,\ell$ and
$0\le n < n_k$ we get, respectively:
\begin{description}
\item[inverse hyperbolic time] 
  $S_{n_k-n}\phi^{cs}(f^nx) \le -7c_s(n_k-n)/8$;
\item[reverse hyperbolic time]
  $S_{n_k-n}\phi^{cs}(f^{m_i-n_k}x)\le -7c_s(n_k-n)/8$.
\end{description}
In the first case we have for inverse $\sigma_s$-hyperbolic times with
$\sigma_s:=e^{-7c_s/8}$
\begin{align*}
  \|Df^{n_k-n}\mid E^{cs}_{f^{n_k-n+1}x}\|
  \le
  \prod\nolimits_{j=n+1}^{n_k}\|Df\mid E^{cs}_{f^j x}\|
  \le
  e^{-7c_s(n_k-n)/8}=\sigma_s^{n-n_k},
\end{align*}
which were implicitly used in~\cite[Proposition 6.4]{ABV00}.  In the
second case we have
\begin{align*}
  \|Df^{n_k-n}\mid E^{cs}_{f^{m_j-n_k}x}\|
  \le
  \prod\nolimits_{j=n+1}^{n_k}
  \|Df\mid E^{cs}_{f^{ m_i-j }x}\|
  \le
  e^{-7c_s(n_k-n)/4}=\sigma_s^{n-n_k}.
\end{align*}
The iterates $m_i-n_k$ are \emph{reverse hyperbolic times} for the
$f$-orbit of $x$ with respect to $m_i$; similar times were used in
\cite{Man88} by Ma\~n\'e and by Liao in~\cite{Li80}. 
  
Pliss' Lemma ensures that there are infinitely many inverse/reverse
hyperbolic times $n_i$ along the $f$-orbit of $x$ with respect to
$m_i$ and, because $\theta_s>0$, we can assume that
$(m_i-n_i)\nearrow\infty$.

\begin{remark}[Chaining property of reverse hyperbolic
  times]\label{rmk:chaining}
  We note that if $n_k$ is a reverse hyperbolic time with respect to
  $m_i$, then it is also a reverse hyperbolic time with respect to all
  times $m$ strictly between $n_i$ and $m_i$ ($n_i<m<m_i$).

  Moreover, if $n_i$ is a reverse hyperbolic time with respect to
  $m_i$ and $n_i<n_j<m_i$ is a reverse hyperbolic time with respect to
  $m_{i+1}>m_i$, then $n_i$ becomes a reverse hyperbolic time with
  respect to $m_{i+1}$.
\end{remark}

Thus, if $h$ is a reverse hyperbolic time with respect to $m_i$, then
$ \|Df^{j}\mid E^{cs}_{f^hx}\|\le\sigma_s^{j} $ for all
$j=1,\dots,m_i-h$ which, roughly speaking, is a hyperbolic time in the
reverse time direction.  This uniform contractive property can be
extended to a neighborhood of the orbit along the center-stable
direction following the same arguments of the proofs of the previous
results for $\sigma_u$-hyperbolic times by replacing backward
contraction with forward contraction; see e.g.~\cite{ABV00}
and~\cite[Lemma 2.2]{Ara2020} and the lower half of
Figure~\ref{fig:revhyptime}.

\begin{proposition}[Pre-disks at reverse hyperbolic times with
  positive frequency]
  \label{pr:reverhyp}
  There exists $\theta_s\in(0,1]$ and $n_0>1$ such that for every
  $x \in U$ and $n> n_0$ with $S_n\phi^{cs}(x)<-c_s n$, there exist
  $l\ge \theta_s\cdot n$ reverse $\sigma$-hyperbolic times
  $1 \le n_1 < \cdots < n_l \le n$ for \( x \) with respect to $n$,
  where $\sigma=e^{-7c_s/8}$.  Moreover, for $\Delta\subset U$ a
  $cs$-disk of radius $\delta_1$ around $f^nx$ and each $i=1,\dots,l$,
  there exists a neighborhood $V_n$ of $f^nx$ in $\Delta$ such that
  \begin{enumerate}
  \item $f^{-(n-n_i)}$ maps $V_n$ diffeomorphically onto a $cs$-disk
    $\Delta_{n_i}=f^{-(n-n_i)}V_n$ of radius $\delta_1$ around
    $f^{n_i}x$;
  \item for every $1\le k \le n-n_i$ and $y, z\in \Delta_{n_i}$,
    \begin{align*}
      \dist_{f^{n_i-n+k}(V_n)}\big(f^k(y),f^k(z)\big)
      \le
      (\sigma^{1/2})^k  \dist_{\Delta_{n_i}}(y,z).
    \end{align*}
  \end{enumerate}
\end{proposition}

\begin{remark}[No pre-disks at inverse hyperbolic times]
  \label{rmk:bad!}
  The same reasoning to construct pre-disks at (reverse) hyperbolic
  times \emph{does not apply to inverse hyperbolic times}, since we
  might have to shrink the domain of the contractions as we move
  backward, so that $cs$-disk centered at $x_{n-k}$ might have a
  radius much smaller than $\delta_1$; see the upper part of
  Figure~\ref{fig:revhyptime}.
\end{remark}

\begin{remark}[Simultaneous hyperbolic times]
  \label{rmk:simultaneous}
  For a possibly smaller neighborhood $\cV$ in the statement of
  Theorem~\ref{thm:NUEnotPH}, it can be show that we have simultaneous
  hyperbolic times and inverse/reverse hyperbolic with positive
  frequency $\theta_u+\theta_s-1$ for all $g\in\cV$ and
  $\leb$-a.e. $x\in M$; see e.g.~\cite[Proposition 6.5]{ABV00}. We
  generalize this idea to intersection of \emph{coherent blocks} in
  the proof of Theorem~\ref{mthm:PesinC1+} in
  Section~\ref{sec:coherent-blocks-hype}.
\end{remark}

\begin{remark}[Roughness of hyperbolic times]
  \label{rmk:delta1simult}
  If $\delta_1>0$ satisfies~\eqref{eq:delta1} for
  $\sigma=\sigma_0\in(0,1)$, then~\eqref{eq:delta1} also holds for all
  $\sigma\in(\sigma_0,1)$.  In what follows we assume, without loss of
  generality, that $\delta_1>0$ is chosen so that~\eqref{eq:delta1}
  holds simultaneously for $\sigma=\sigma_s$ and
  $\sigma=\sigma_u$.
  % Moreover, for future use, we also assume that
  % $\delta_1^{\zeta}(1-\sigma_s^{\zeta/2})^{-1}<c_u/4$.
\end{remark}

\subsection{Schedules and coherent blocks}% and synchronization}
\label{sec:schedul-omega-limits}

The following results from~\cite{Pinho2011} and~\cite{pinheiro20} will
be used as tools in the proofs of the main theorems.

A \emph{schedule} of events is a measurable map $\cU:U\to 2^{\ZZ^+_0}$
which is \emph{asymptotically invariant} if for $x\in U$
\begin{enumerate}
\item $\#\cU(x)=\infty$; and
\item $\cU(x)\cap[n,+\infty)=\cU(f(x))\cap[n,+\infty)$ for
  every big $n\in\ZZ^+$.
\end{enumerate}
The asymptotically invariant schedule $\cU=(\cU(x))_{x\in U}$ has {\em
  positive frequency} if for each $x\in U$ it satisfies
\begin{align*}
  d^+(\cU(x))
  :=
  \limsup_{n\nearrow\infty}\frac1n\#\big(\cU(x)\cap[0,n)\big)
  >0.
\end{align*}
A schedule of events $\cU=(\cU(x))_{x\in U}$ is \emph{coherent} if it
satisfies the following properties:
\begin{enumerate}
\item if $n\in\cU(x)$ then $n-j\in\cU(f^j(x))$ for every $x\in\U$ and
  $n>j\ge0$; and
\item if $n\in\cU(x)$ and $m\in\cU(f^n(x))$, then $n+m\in\cU(x)$ for
  every $x\in U$ and $n,m\ge1$.
\end{enumerate}

\begin{remark}\label{rmk:excoherent}
  The schedules of events $\cU_1,\cU_2:M\to2^{\ZZ^+_0}$ given by,
  respectively, inverse hyperbolic times and hyperbolic times, are all
  $f$-coherent schedule of events with positive frequency.
\end{remark}

We define the \emph{$f$-coherent block for $\cU$} or, for short, the
\emph{$\cU$-block}, as
\begin{align*}
  B_{\cU}
  =
  \big\{x\in\cap_{n\ge0}f^{n}(U):
  j\in\cU(f^{-j}(x)), \,\forall j\ge0\big\}.
\end{align*}

\begin{theorem}\cite[Theorem 6.4]{pinheiro20}\label{thm:blockpositive}
  If $\mu$ is an ergodic $f$-invariant probability on $U$ and
  $\U:U\to 2^{\ZZ_0^+}$ is a coherent schedule, then
  $\mu(B_{\cU}) = d^+(\cU(x))$ for $\mu$-almost every $x\in U$.
\end{theorem}

\subsection{Coherent block for reverse hyperbolic times}
%{Long reverse hyperbolic times and stable leaves of uniform size}
\label{sec:long-reverse-hyperb}

\begin{figure}[htpb]
  \centering \includegraphics[width=9cm]{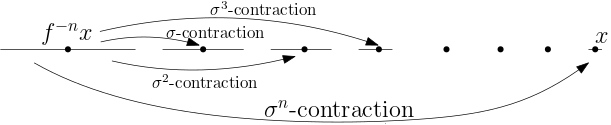}
  \caption{\label{fig:revhyptime-1}Forward contractions from $f^{-n}x$
  to $x$.}
\end{figure}

We note that since we have a physical $f$-invariant ergodic probability
measure $\mu$, then the limit~\eqref{eq:NUC} holds for $\mu$-a.e. $x$
also for the inverse transformation $f^{-1}$, that is
\begin{align*}
  \lim_{n\to+\infty} \frac1n\sum\nolimits_{j=1}^n\phi^{cs}(f^{-j}x)<-c_s,
  \quad \mu-\text{a.e.} x.
\end{align*}
We may then find ``hyperbolic times'' in this setting, that is,
times $n\ge1$ so that 
\begin{align*}
  \sum\nolimits_{i=0}^{k-1} \phi^{cs}(f^{-(n-k+i)}x) < -7 k c_s/8, \quad 0<k\le n;
\end{align*}
or equivalently
\begin{align*}
  \|Df^k\mid E^{cs}_{f^{-n}x}\|
  \le
  \prod\nolimits_{i=0}^{k-1}\|Df\mid E^{cs}_{f^{-n+k-i}x}\| \le e^{-k \cdot
  7c_s/8},
  \quad 0< k \le n.
\end{align*}
This means that  $-n$ \emph{becomes a reverse hyperbolic time with
  respect to $0$}; see Figure~\ref{fig:revhyptime-1}.

For $\mu$-a.e $x$ the family of absolute values of such times can be
seen as a schedule $\wh{\cU}$ with respect to the dynamics of
$g:=f^{-1}$ which is coherent and has positive frequency. Points $x$
in the corresponding reverse hyperbolic block $B_{\wh{\cU}}$ satisfy
$j\in\wh{\cU}(g^{-j}x)$ for all $j\ge0$. That is, $x\in B_{\wh{\cU}}$
if, and only if, $j\in \wh{\cU}(f^jx)$ and so $0$ \emph{becomes a
  reverse hyperbolic time with respect to all $j>0$}.

\begin{figure}[htpb]
  \centering \includegraphics[width=9cm]{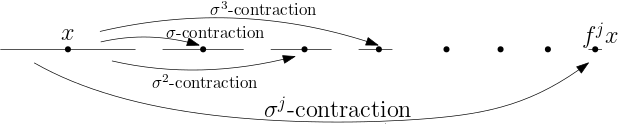}
  \caption{\label{fig:revhyptimelong}Forward contractions from $x$
  to $f^jx$ for any $j>0$.}
\end{figure}

We say that $x\in B_{\wh{\cU}}$ is a point with a \emph{long reverse
  hyperbolic time}; see Figure~\ref{fig:revhyptimelong}. 

\subsubsection{Stable leaves of uniform size}
\label{sec:stable-leaves-unifor}

% Let us fix a $cu$-disk $\Sigma$ that is non-uniformly expanding along
% the $E^{cu}$ direction and non-uniformly contracting along the
% $E^{cs}$ direction on a subset $H\subset\Sigma$ with full
% $\leb_\Sigma$-measure.

We assume, without loss of generality, that $\delta_1>0$
from~\eqref{eq:delta1} is such that the exponential map of $M$ is
invertible on $\delta_1$-balls in the tangent space, i.e.
$\exp_x: B(0,\delta_1)\subset T_xM\to M$ is a diffeomorphism with its
image, where $B(0,\delta_1)$ denotes $\{w\in T_xM: \|w\|<\delta_1\}$
and, additionally, that $\exp_x(E^{cs}_x\cap B(0,\delta_1/2))$ is a
$cs$-disk for any given $x\in M$.

For each $x\in B^s:=B_{\wh{\cU}}$ we have a reverse
$\sigma_s$-hyperbolic time for $f^ix$, for each $i=1,2,\ldots$

For each $i>1$, choosing a $cs$-disk
$\Delta_i=\exp_{f^ix}\big(E^{cs}_{f^ix}\cap
B(0,\delta_1/2)\big)\subset U$ at $f^ix$ we get, from
Proposition~\ref{pr:reverhyp}, a neighborhood $V_i$ of $f^ix$
in $\Delta_i$ so that for each $k=1,\dots,i$
\begin{itemize}
\item $f^{-k}V_{\ell_i}$ is a $cs$-disk through $f^{i-k}x$ with radius
  at most $\delta_1$; and
\item $f^k:f^{-k}V_i\to\Delta_i$ is a $\sigma_s^{k/2}$-contraction.
\end{itemize}
If we set $D_i:=f^{-i}V_i$ then, by the Ascoli-Arzela Theorem, there
exists a $cs$-disk $D_x$ with radius $\delta_1$ around $x$,
which is an accumulation point of $D_i$ in the $C^1$ topology when
$i\nearrow \infty$.

Moreover, by continuity of the map $f$, we have that $f^k\mid D_x$ is
a diffeomorphism from the $cs$-disk $D_x$ into the $cs$-disk
$D_x^k:=f^k D_x$, and
$\big(f^k\mid D_x\big)^{-1}=f^{-k}\mid D_x^k:D_x^k\to D_x$ is a
$\sigma_s^{-k/2}$-expansion for each $k\ge1$.

It follows from, e.g. the Non-Uniform Hyperbolic Theory for hyperbolic
measures with dominated splitting \cite[Proposition 8.9]{Abdenur2011},
that $D_x$ is the stable manifold at $x$ with radius $\delta_1$, that
is, $D_x=W^s_x(\delta_1)$ and $T_yD_x=E^{cs}_y$ for all $y\in D_x$. We
have proved the following.

\begin{proposition}[Long stable leaves on the reverse hyperbolic block
  $B_{\wh{\U}}$]
  \label{pr:longstable}
  For $x\in B_{\wh{\U}}$ there exist a center-stable disk
  $W^s_x(\delta_1)$ tangent to the center-stable direction and with
  radius $\delta_1>0$ centered at $x$, which is the local stable
  manifold. More precisely
  $ W^s_x(\delta_1) = \{y\in M:
  d(f^ky,f^kx)\le\delta_1\sigma_s^{k/2}, k\ge0\}.$
\end{proposition}

\subsection{Unstable leaves of uniform size}
\label{sec:unstable-leaves-unif}

Analogously, we obtain local unstable manifolds throught every point
of the coherent block $B^u:=B_{\U_u}$, given by the coherent schedule $\U_u$
of $\sigma_u$ hyperbolic times defined at $\mu$-a.e. point $x$, where
$\mu$ is a physical/SRB measure for $f$.

\begin{proposition}[Long unstable leaves on the hyperbolic block
  $B^u$]
  \label{pr:unstable}
  For $x\in B^u$ there exist a center-unstable disk
  $W^u_x(\delta_1)$ tangent to the center-unstable direction and with
  radius $\delta_1>0$ centered at $x$, which is the local unstable
  manifold. More precisely
  $ W^u_x(\delta_1) = \{y\in M:
  d(f^ky,f^kx)\le\delta_1\sigma_u^{k/2}, k\le0\}.$
\end{proposition}

\begin{proof}
  Each point $x\in B^u$ is such that every $n\in\ZZ^+$ is a
  $\sigma_u$-hyperbolic time for $f^{-n}(x)$. We can apply
  Lemma~\ref{le:contraction} starting at a $cu$-disk $\Delta$ at
  $f^{-n}(x)$ to obtain a $cu$-disk $\Delta_n$ of radius $\delta_1>0$
  around $x$ which is uniformly contracted backwards at a rate
  $\sigma_u^{1/2}$ for up to $n$ iterates. Just like in the proof of
  Proposition~\ref{pr:longstable}, we conclude that these disks $D_n$
  accumulate to an unstable disk $W^{cu}_x(\delta_1)$; see also
  e.g. \cite[Lemma 3.7]{ABV00}.
\end{proof}

%%%%%%%%%%%%%%%%%%%%%%%%%%%%%%%%%%%%%%%%%%%%%%%%%%%%%%

\section{Long (un)stable leaves on subsets with positive measure}
\label{sec:long-unstable-leaves}

Here we prove Theorem~\ref{mthm:PesinC1+}. We start by showing that
each hyperbolic ergodic measure is automatically a non-uniformly
hyperbolic measure in the sense of
Subsection~\ref{sec:non-uniform-expans} for a power of the map. This
holds for any ergodic hyperbolic and dominated invariant probability
measure for a $C^1$ diffeomorphism.

We then use this in the $C^{1+}$ setting to take advantage of the
existence of hyperbolic times in different versions to construct
(un)stable manifolds with uniformly bounded size on coherent blocks
with positive measure.

\subsection{Hyperbolic dominated measures and
  non-uniform hyperbolic dominated splitting}
\label{sec:hyperb-dominat-measu-1}

The following shows that hyperbolic and dominated measures have a
non-uniformly hyperbolic splitting for a power of the dynamics.
%\margem{enough to use continuous splitting here ?}
\begin{lemma}[Non-uniform contraction for a power]{\cite[Lemma
    8.4]{Abdenur2011}}
  \label{le:NUCpower}
  Let $f$ be a $C^1$ diffeomorphism, $\mu$ be an ergodic $f$-invariant
  probability measure and $E^{cs}\subset T_{\supp \mu}M$ be a
  $Df$-invariant continuous subbundle defined over $\supp\mu$. Let
  $\lambda^+_{cs}$ be the upper Lyapunov exponent in $E^{cs}$ of the
  measure $\mu$ as in~\eqref{eq:middlexp}.  Then, for any
  $\epsilon > 0$ with $\lambda^+_{cs}+\epsilon<0$, there exists an
  integer $N(\epsilon,\mu)$ such that, for $\mu$-a.e. $x$ and each
  $N\ge N(\epsilon,\mu)$, the Birkhoff averages
  $S^{f^N}_k\phi_N^{cs}(x) / Nk$ converge towards a number contained
  in $[\lambda^+_{cs},\lambda^+_{cs}+\epsilon)$ when
  $k\nearrow\infty$.
\end{lemma}

\begin{remark}[Non-uniform contraction for a power and dependence of
  $\epsilon$]
  \label{rmk:NUCpower}
  Therefore, if $\lambda_{cs}^+<0$, then $E^{cs}$ becomes non-uniform
  contracting $\mu$-a.e. for a power $f^N$, where $N=N(\epsilon,\mu)$
  and $\epsilon>0$ so that $\lambda_{cs}^++\epsilon<0$.  The proof
  of~\cite[Lemma 8.4]{Abdenur2011} shows that $N=N(\epsilon,\mu)$ is
  determined by the condition
  $\mu(\phi^{cs}_N)<N(\lambda^+_{cu}+\epsilon)<0$ and so 
  $N(\epsilon,\mu)\nearrow\infty$ when $\epsilon\searrow0$ (following
  Kingman's Subadditive Ergodic Theorem~\cite[Section 3.3]{ViOl16}).
\end{remark}

Recalling~\eqref{eq:middlexp}, if $\lambda_{cu}^-<0$, then % by the
% invertibility of $f$ and Oseledets' Multiplicative Ergodic Theorem we
% get
% $\lim_{n\nearrow\infty} \log\|Df^{-n}\mid
% E^{cu}_x\|^{1/n}=\lambda_{cu}^-<0$. Therefore,
replacing %$f$ by $f^{-1}$ and
$\phi_N^{cs}(x)$ by $\phi_N^{cu}(x)$ in the statement of
Lemma~\ref{le:NUCpower}, we conclude that for any $\epsilon>0$ there
exists $\wt{N}(\epsilon,\mu)\in\ZZ^+$ such that for $\mu$-a.e $x$ and
each $N\ge \wt{N}(\epsilon,\mu)$, the averages
$S_k^{f^{-N}}\phi^{cu}_N(x)/Nk$ converge towards a number in
$[\lambda^-_{cu},\lambda^-_{cu}+\epsilon)$.

Hence, if $\lambda_{cu}^-<0$, then $E^{cu}$ becomes non-uniform
expanding $\mu$-a.e. with respect to a power $f^{N}$, where
$N=\wt{N}(\epsilon,\mu)$ for $\epsilon>0$ so that
$\lambda_{cu}^-+\epsilon<0$.  Altogether, the threshold $N$ ultimately
depends on $\mu$, $|\lambda_{cs}^+|$ and $|\lambda_{cu}^{-}|$, and we
obtain the following.

\begin{proposition}[Hyperbolic dominated measure is non-uniformly
  hyperbolic]
  \label{pr:NUEfrompos}
  Let $f$ be a $C^1$ diffeomorphism, $\mu$ be an ergodic $f$-invariant
  probability measure and $T_{\supp \mu}M=E^{cs}\oplus E^{cu}$ be a
  $Df$-invariant and dominated splitting over $\supp\mu$ such that
  $\max\{\lambda_{cs}^+, \lambda_{cu}^-\}<0$. Then there exists
  $N=N(f,\mu,|\lambda_{cs}^+|,|\lambda_{cu}^{-}|)\in\ZZ^+$ so that
  $f^N$ is non-uniformly hyperbolic with respect to $\mu$, that is
  both~\eqref{eq:NUE0} and~\eqref{eq:NUC} hold on a full $\mu$-measure
  subset with respect to the iterates of $g:=f^N$.
\end{proposition}

\subsection{Coherent blocks and hyperbolic times}
\label{sec:coherent-blocks-hype}

We are now ready to prove Theorem~\ref{mthm:PesinC1+}.

\begin{proof}[Proof of Theorem~\ref{mthm:PesinC1+}]
  At this point we have $g=f^N$ which is non-uniformly hyperbolic on a
  full $\mu$-measure subset. However, $\mu$ might not be
  $g$-ergodic. From \cite[Lemma 3.13]{Pinho2011} since $\mu$ is
  $f$-ergodic we decompose
  \begin{align}
    \label{eq:ergpower}
    \mu=\frac1k\big(\nu+f_*\nu+\cdots+f_*^{k-1}\nu\big),
  \end{align}
  where $k\in\ZZ^+$
  divides $N$ and $\nu$ is $f^k$-invariant and $g$-ergodic.

  Hence, using the asymptotically invariant and coherent schedules
  $\wh{\U}$ of reverse $\sigma_s$-hyperbolic times for $g^{-1}$
  (recall Subsection~\ref{sec:long-reverse-hyperb}) and $\U_u$ of
  $\sigma_u$ hyperbolic times, defined for $\nu$-a.e.
  $x\in\supp(\mu)$, we obtain from Theorem~\ref{thm:blockpositive}
  that the corresponding $g$-coherent blocks $B^s:=B_{\wh{\U}}$ and
  $B^u:=B_{\U^u}$ satisfy
  \begin{align*}
    \nu(B^s)=d^+(\wh{\U})\ge\theta_s
    \qand
    \nu(B^u)=d^+(\U_u)\ge\theta_u.
  \end{align*}
  Here $\theta_s,\theta_u\in(0,1)$ are the lower bounds for the
  asymptotic density of Pliss times, which depend on $g$ and the
  values $\sigma_s=\exp(\lambda^+_{cs}+\epsilon)^{7/8}$ and
  $\sigma_u=\exp(-\lambda_{cu}^-+\epsilon)^{7/8}$ from the proof of
  Propostion~\ref{pr:NUEfrompos}. More precisely, we have
  \begin{align*}
    \theta_u=\theta_u(Df,N,|\lambda_{cu}^-|)
    &=
      \frac{|\log\sigma_u|}{8\log\sup_{x\in U}\|(Df^N\mid
    E^{cu}_x)^{-1}\|^{-1}-7|\log\sigma_u|}
      \qand
    \\
    \theta_s=\theta_s(Df,N,|\lambda_{cs}^+|)
    &=
      \frac{|\log\sigma_s|}{8\log\sup_{x\in U}\|Df^N\mid
    E^{cs}_x\|^{-1}-7|\log\sigma_s|}.
  \end{align*}
  On the one hand, item (1) of the statement of
  Theorem~\ref{mthm:PesinC1+} follows from
  Proposition~\ref{pr:longstable}, where the inner radius of
  $W^s_x(\delta_1)$ for each $x\in B^s$ is
  $\delta_1=\delta_1(f,N,|\lambda_{cs}^+|)$, since we have long
  reverse $\sigma_s$-hyperbolic times by definition of coherent block;
  recall Subsection~\ref{sec:long-reverse-hyperb}.
  
  On the other hand, from Proposition~\ref{pr:unstable}, we have
  uniformly sized unstable manifolds $W^u_y(\delta_1)$ through each
  point $y$ of $B^u$, where $\delta_1=\delta_1(f,N,|\lambda_{cu}^-|)$.
  This proves item (2) of the statement of
  Theorem~\ref{mthm:PesinC1+}.
  
  Since $\nu$ is $g$-ergodic, there exists $\ell\in\ZZ_0^+$ so that
  $\nu(B^u\cap f^{-\ell}B^s)>0$. Setting $B$ as in item (3) of the
  statement of Theorem~\ref{mthm:PesinC1+}, we complete the proof by
  noting that each $x\in\cH$ reaches $B^s$ in at most $\ell$ iterates,
  and so there are constants $c,C>0$ as stated.

  Finally, for the regularity of the lamination $\cF^s$, the absolute
  continuity and H\"older continuity of the Jacobian of holonomy maps
  follow in general as in~\cite[Chapter 8, Theorems 8.6.1 \&
  8.6.15]{BarPes2007}.

  More precisely, with our stronger assumptions, we have that for each
  pair of non-intersecting $cu$-disks $\gamma_1,\gamma_2$ crossing
  $\cF^s_z$ (necessarily transversely and with angles bounded away
  from zero, as a consequence of the dominated splitting) then, after
  setting $F^s:=\cup_{x\in B^s}\cF^s_x=\cup_{x\in B^s}W^s_x(\delta_1)$
  and the holonomy $\Theta:\gamma_1\cap F^s\to\gamma_2\cap F^s$ given
  by $\Theta(x):=\cF^s_x\cap\gamma_2$, we have
  $\Theta_*\leb_{\gamma_1\cap F^s}\ll\leb_{\gamma_2}$. Moreover, the
  corresponding density
  $\rho=\rho_{\gamma_1,\gamma_2}=\frac{d\big(\Theta_*\leb_{\gamma_1\cap
      F^s}\big)}{\leb_{\gamma_2}}$ is given by
  \begin{align}\label{eq:defrho}
    \rho(x):=
    \exp\sum\nolimits_{i\ge0}\big(J^{cu}(f^i\Theta(x))-J^{cu}(f^ix)\big).
  \end{align}
  Since $x$ and $\Theta(x)$ belong to the same local stable leaf
  $W^s_z(\delta_1)$ for some $z\in B^s$, and $J^{cu}$ is a
  $\eta$-H\"older for some $\eta\in(0,1]$, we can find a constant
  $C_J>0$ so that
    \begin{align}
      &\sum\nolimits_{i\ge0}\big|J^{cu}(f^i\Theta(x))-J^{cu}(f^ix)\big|
        \le
        C_J\sum\nolimits_{i\ge0}
        \dist_{\cF^s_{f^ix}}\big(f^i\Theta(x),f^ix\big)^\eta
        \nonumber
      \\
      &\le 
        C_J\dist_{\cF^s_x}(x,\Theta(x))^\eta
        \sum\nolimits_{i\ge0}\sigma_s^{2i\eta/3}
        \le
        C_J\dist_{\cF^s_x}(x,\Theta(x))^\eta/(1-\sigma_s^{2\eta/3}).
        \label{eq:bddrho}
    \end{align}
    Since $\dist_{\cF^s_x}(x,\Theta(x))\le\delta_1$ and is bounded
    away from zero for all $x\in\gamma_1$, we conclude that $\rho(x)$
    is bounded above and below away from zero and infinity, as stated.
  \end{proof}

  % \begin{proof}[Proof of Corollary~\ref{mcor:hWus}]
  %   In the setting of Theorem~\ref{mthm:PesinC1+}, any given
  %   $\mu\in DH(f)$ with $Df$-invariant dominated splitting
  %   $T_{\supp(\mu)}M=E^{cs}\oplus E^{cu}$ and $d_{cu}=\dim E^{cu}>0$,
  %   $d_{cs}=\dim E^{cs}>0$ (since $\mu$ is non-atomic) satisfies, from
  %   the Margulis-Ruelle Inequality~\cite[Section 4.12]{Man87}
  %   \begin{align*}
  %     h_\mu(f)\le d_{cu}\cdot|\lambda_{cu}^-|
  %     \qand
  %     h_\mu(f)\le d_{cs}\cdot|\lambda_{cs}^+|.
  %   \end{align*}
  %   Hence, if $h_\mu(f)>h$ and $d=\dim M$, then
  %   $ |\lambda_{cu}^-|>\frac{h}{d_{cu}}\ge\frac{h}{d-1} \qand
  %   |\lambda_{cs}^+| >\frac{h}{d_{cs}}\ge\frac{h}{d-1}.$ Thus, from
  %   the proof of Theorem~\ref{mthm:PesinC1+},
  %   Remark~\ref{rmk:NUCpower} and Proposition~\ref{pr:NUEfrompos}, we
  %   find a lower bound $\delta_1$ for the inner radius of the
  %   (un)stable manifolds on the point of coherent blocks $B^s$ and
  %   $B^u$ so that whenever $\delta<\delta_1$
  %   \begin{align*}
  %     L^s(\delta)
  %     \supset
  %     B^s \qand \mu(B^s)>\theta_s;
  %     \qquad
  %     L^u(\delta)
  %     \supset
  %     B^u \qand \mu(B^u)>\theta_u;
  %   \end{align*}
  %   where $\delta_1,\theta_s,\theta_u$
  %   depend only on $f$, the largest contracting Lyapunov exponent
  %   $\lambda_{cs}^+$ and the smallest expanding Lyapunov exponent
  %   $\lambda_{cu}^-$, whose modulus are both bounded below by
  %   $h/(d-1)$.  The statement of the corollary follows.
  % \end{proof}

%%%%%%%%%%%%%%%%%%%%%%%%%%%%%%%%%%%%%%%%%%%%%%%%%%%%%%%%%%%%%%%%%

\section{GMY structure for non-uniformly hyperbolic attracting sets}
\label{sec:existence-physic-mea-1}

Here we prove Theorem~\ref{mthm:NUHypGMY} and
Corollaries~\ref{mcor:abv} and~\ref{mcor:GMYexp} following the same
strategy presented in~\cite[Chapter 7]{Alves2020b} and also used
in~\cite{ADLP,AlPi10}, citing and adapting the main tools according to
our more general assumptions.

We start by recalling the definion of a GMY structure, in
Subsection~\ref{sec:gibbsm-struct}. In
Subsection~\ref{sec:constr-gmy-struct-1}, we describe how to obtain
this structure in our dynamical setting, preparing the proof of
Theorem~\ref{mthm:NUHypGMY}
% and Corollary~\ref{mcor:abv}.
by constructing the family of unstable disks in a cylinder.  In
Subsection~\ref{sec:recurr-cylind-over}, we use syncronization and the
stable coherent block to build the family of stable disks in the same
cylinder obtained in Subsection~\ref{sec:constr-gmy-struct-1}.  In
Subsection~\ref{sec:existence-full-gmy}, we prove
Theorem~\ref{mthm:NUHypGMY} and Corollary~\ref{mcor:abv}.

\subsection{Gibbs-Markov-Young structure}
\label{sec:gibbsm-struct}

We give here the precise definitions combining recent developments
from~\cite{AlPi10,AlLi15} and~\cite{Alves2020b}.

If $u=\dim E^{cu}$ and $s=\dim E^{cs}$ we write $D^s, D^u$ for the
unit compact balls on $\RR^s$ and $\RR^u$, respectively, and say that
any diffeomorphic image of $D^u\times D^s$ is a \emph{cylinder}. 

We say that $\Gamma^u=\{\gamma^u\}$ is a \emph{continuous family of
  $C^1$ unstable manifolds} if there is a compact set~$K^s$, a unit
disk $D^u$ of some $\RR^n$, and a map
$\Phi^u\colon K^s\times D^u\to M$ such that
\begin{enumerate}[(i)]
\item $\gamma^u=\Phi^u(\{x\}\times D^u)$ is an unstable
manifold;
\item $\Phi^u$ maps $K^s\times D^u$ homeomorphically onto its
  image;
\item $x\mapsto \Phi^u\vert(\{x\}\times D^u)$ defines a
continuous map from $K^s$ into $\text{Emb}^1(D^u,M)$.
\end{enumerate}
Here $\text{Emb}^1(D^u,M)$ denotes the space of $C^1$ embeddings from
$D^u$ into $M$.  Continuous families of $C^1$ stable manifolds are
defined similarly.

We say that  a set $\Lambda\subset M$ has a \emph{hyperbolic product
structure} if there exist a continuous family of local unstable manifolds
$\Gamma^u=\{\gamma^u\}$ and a continuous family of local  stable manifolds
$\Gamma^s=\{\gamma^s\}$ such that
\begin{enumerate}[(1)]
    \item $\Lambda=(\cup \gamma^u)\cap(\cup\gamma^s)$;
    \item $\dim \gamma^u+\dim \gamma^s=\dim M$;
    \item each $\gamma^s$ intersects each $\gamma^u$ in exactly
      one point;
    \item stable and unstable manifolds are transversal with angles
      bounded away from~$0$.
  \end{enumerate}
If $\Lambda\subset M$ has a product structure, we say that
$\Lambda_0\subset \Lambda$ is an {\em $s$-subset} if $\Lambda_0$ also
has a product structure and its defining families $\Gamma_0^s$ and
$\Gamma_0^u$ can be chosen with $\Gamma_0^s\subset\Gamma^s$ and
$\Gamma_0^u=\Gamma^u$; {\em $u$-subsets} are defined analogously. For
convenience we shall use the following notation: given $x\in\Lambda$,
let $\gamma^{*}(x)$ denote the element of $\Gamma^{*}$ containing $x$,
for $*=s,u$. Also, for each $n\ge 1$ let $(f^n)^u$ denote the
restriction of the map $f^n$ to $\gamma^u$-disks and let
$\det D(f^n)^u$ be the Jacobian of $D(f^n)^u$.

We say that \( f \) admits a \emph{Gibbs-Markov-Young (GMY) structure}
if there exist a set \( \Lambda \) with hyperbolic product structure
satisfying the following additional properties.
\begin{enumerate}[(I)]
\item \emph{Detectable}: $\leb_\gamma(\Lambda)>0$ for each
  $\gamma\in\Gamma^u$.
\item \emph{Markov}: there are pairwise disjoint $s$-subsets
  $\Lambda_1,\Lambda_2,\dots\subset\Lambda$ such that
   \begin{enumerate}
   \item
     $\leb_{\gamma}\big((\Lambda\setminus\cup\Lambda_i)\cap\gamma\big)=0$
     on each $\gamma\in\Gamma^u$.
   \item for each $i\in\NN$ there is $R_i\in\NN$ such that
     $f^{R_i}(\Lambda_i)$ is $u$-subset, and for all $x\in \Lambda_i$
         $$
         f^{R_i}(\gamma^s(x))\subset \gamma^s(f^{R_i}(x))\qand
         f^{R_i}(\gamma^u(x))\supset \gamma^u(f^{R_i}(x)).
         $$
       \end{enumerate}
 \end{enumerate}
 The Markov property enables the definition of a recurrence time
 $R:\Lambda\to\ZZ^+$ and return map $f^R:\Lambda\to\Lambda$ defined on
 a full $\leb_\gamma$-measure subset $\Lambda\cap\gamma$ for each
 $\gamma\in\Gamma^u$ so that
\begin{align*}
  R\mid\Lambda_i\equiv R_i
  \qand
  f^R\mid\Lambda_i \equiv f^{R_i}\mid \Lambda_i.
\end{align*}
Hence, there is a subset $\Lambda'\subset\Lambda$ intersecting each
$\gamma\in\Gamma^u$ in a full $\leb_\gamma$-measure subset of
$\gamma\cap\Lambda$ such that $(f^R)^n(x)$ lies in some $\Lambda_i$
for each $n\ge0$ and all $x\in\Lambda'$. For $x,y\in\Lambda'$ we set
the \emph{separation time}\footnote{We convention that
  $\min\emptyset=\infty$ and set $s(x,y)=0$ for points in
  $\Lambda\setminus\Lambda'$.}
$s(x,y) := \min\{n\ge0: (f^R)^n(x) \;\&\; (f^R)^n(y) \text{ belong to
  different } \Lambda_i\}$.
  
The next conditions assume that there are constants $C>0$ and
$0<\beta<1$, depending on $f$ and $\Lambda$, satisfying the following.
\begin{enumerate}[(I),resume]
  \item \emph{Contraction on stable leaves}: for all 
    $\gamma^s\in\Gamma^s$, $x,y\in\gamma^s\cap\Lambda_i$ 
    \begin{enumerate}[(a)]
    \item $\dist\big((f^R)^n(x),(f^R)^n(y)) \le C\beta^n$ for all
      $n\ge0$; and
    \item $\dist(f^n(y),f^n(x)) \le C d(y,x)$ for all $1\le n\le R_i$.
    \end{enumerate}
  \item \emph{Expansion on unstable leaves}: for each $i\ge1$ and all
    $\gamma^u\in\Gamma^u$, $x, y \in \Lambda_i\cap\gamma^u$
    \begin{enumerate}[(a)]
    \item
      $\dist((f^R)^n(y),(f^R)^n(x))\le C\beta^{s(x,y)-n}$ for all $n\ge0$; and
    \item $\dist(f^i(y),f^i(x))\le C\dist(f^R(y),f^R(x))$ for all
      $0< i\le R=R(\Lambda_i)$.
    \end{enumerate}
 \item \emph{Bounded distortion}: for all $i\ge1$,
   $\gamma^u\in\Gamma^u$ and $x, y \in \Lambda_i\cap\gamma^u$
  $$
  \log\frac{\det D(f^{R_i})^u(x)}{\det D(f^{R_i})^u(y)}\leq C
  \beta^{s(f^{R}(x),f^{R}(y))}.
$$

\item \emph{Regularity of the stable holonomy}:
  for all $\gamma,\gamma'\in\Gamma^u$ we define
  $\Theta:\gamma\cap\Lambda\to\gamma'\cap\Lambda$ by setting
  $\Theta(x)$ equal to $\gamma^s(x)\cap \gamma'$, and
  $\Theta_*\leb_\gamma$ is absolutely continuous with respect to
  $\leb_{\gamma'}$ and its density $\rho=\rho_{\gamma,\gamma'}$
  satisfies
  \begin{align*}
    \frac1C
    \le
    \int_{\gamma'\cap\Lambda} \rho\,d\leb_{\gamma'}
    \le
    C
    \qand
    \log\frac{\rho(x)}{\rho(y)}\le C\beta^{s(x,y)},
    x,y\in\gamma'\cap\Lambda.
  \end{align*}
        % $$\displaystyle \frac{d(\Theta_*\leb_{\gamma})}{d\leb_{\gamma'}}(x)=
        % \prod_{i=0}^\infty\frac{\det Df^u(f^i(x))}{\det
        % Df^u(f^i(\Theta^{-1}(x)))}.$$
\end{enumerate}

A GMY structure is a \emph{full GMY structure} if every disk in
$\Gamma^u$ is contained in $\Lambda$.

We define a return time function $R: \Lambda\to\NN$ by
$R|\Lambda_i= R_i$ and we say that the GMY structure has
\emph{integrable return times} if
$ \int_{\gamma\cap\Lambda}R\, d\leb_\gamma <\infty$ for some
$\gamma \in\Gamma^u$\footnote{Hence, for all $\gamma \in\Gamma^u$ by
  property (VI).}.

%%%%%%%%%%%%%%%%%%%%%%%%%%%%%%%%%%%%%%%%%%%%%%%

\subsection{Construction of the unstable family}
\label{sec:constr-gmy-struct-1}

The first step is provided by the following known result from Alves,
Bonatti and Viana~\cite{ABV00} and Vasquez~\cite{vasquez2006}.

\begin{theorem}[Dominated non-uniform expansion and $cu$-Gibbs
  states]{\cite[Theorem 6.3]{ABV00} \& \cite[Theorem 3.2 \& Corollary
    4.1]{vasquez2006}}
  \label{thm:abvGibbs}
  Let $f$ be a $C^{1+}$ diffeomorphism admitting an attracting compact
  set $A$ with a dominated splitting $T_AM=E^{cs}_A\oplus E^{cu}_A$.
  Assume that $f$ is non-uniformly expanding along the centre-unstable
  direction in the trapping neighborhood $U$ of $A$, i.e., we have
  condition~\eqref{eq:NUE0} on $H\subset H_u$ with $\leb(H)>0$. Then
  \begin{enumerate}[(A)]
  \item $f$ has some ergodic Gibbs $cu$-state $\mu$ supported in
    $\Lambda$;
  \item every ergodic physical/SRB $f$-invariant probability measure
    supported in $U$ is a $cu$-Gibbs state.
  \end{enumerate}
  More precisely, there exists a cylinder $\cC_0$ and a family
  $\Gamma$ of disjoint $cu$-disks contained in $\cC_0$ which are
  graphs over $D^u$, and a ergodic $f$-invariant probability measure
  $\mu$ supported in $\Lambda$, satisfying
  \begin{enumerate}[(a)]
  \item there exist a $cu$-disk $D$ such that $\leb_D(H)>0$, so that
    \begin{enumerate}[(i)]
    \item each disk $\gamma\in\Gamma$ is accumulated by sub-disks of
      radius $\delta_1$ in $f^n(D)$ around points $f^n(x)$ such that
      $n$ is a $\sigma_u$-hyperbolic time for $x\in D\cap H$ with
      $\sigma_u=e^{-7c_u/8}$; consequently
    \item each disk $\gamma\in\Gamma$ is uniformly backward
      contracted:
      $\dist_{f^{k}\gamma}\big(f^{-k}y,f^{-k}z\big) \le
      \sigma_u^{k/2}\dist_\gamma(y,z)$ for all $y,z\in\gamma$ and
      $k\in\ZZ^+$; and
    \item the $d_{cu}=\dim E^{cu}$ larger Lyapunov exponents of $\mu$
      are larger than $\log \sigma_u^{-1/2}=7c_u/16$;
    \end{enumerate}    
  \item $\cC_0$ contains a ball whose radius $r>0$  depends only on $f$;
  \item there exists $\alpha>0$ so that the union
    $\wh{\Gamma}=\cup_{\gamma\in\Gamma}\gamma$ (of the disks in
    $\Gamma$) satisfies $\mu(\wh{\Gamma})\ge\alpha$;
  \item the restriction of $\mu$ to $\wh{\Gamma}$ has absolutely
    continuous conditional measures along the disks in $\Gamma$:
    for every measurable bounded function $\vfi:M\to\RR$ we have
    \begin{align*}
      \int_{\wh{\Gamma}}\vfi\,d\mu
      =
      \int_{\gamma\in\Gamma}\left( \int_{x\in\gamma}\vfi(x)\rho_\gamma(x)
      \,d\leb_\gamma(x)\right) \, d\wh{\mu}(\gamma)
    \end{align*}
    where $\leb_\gamma$ is the induced volume measure on $\gamma$ from
    $\leb$; and $\wh{\mu}=\pi_*\mu$ is the quotient measure, for
    $\pi:\wh{\Gamma}\to\Gamma$ the natural map
    $x\in\wh{\Gamma}\mapsto \gamma_x\in\Gamma$. In addition, the
    densities $\rho_\gamma$ are bounded away from zero and infinity
    depending only on $f$ and $c_u$ (the rate of non-uniform expansion
    from~\eqref{eq:NUE0}) due to the relation
    \begin{align*}
      \frac{\rho_{\gamma}(x)}{\rho_\gamma(y)}
      =
      \prod_{k\ge0}
      \frac{\det(Df^{-1}\mid E^{cu}_{f^{-k}x})}
      {\det(Df^{-1}\mid E^{cu}_{f^{-k}y})}, \quad x,y\in\gamma.
    \end{align*}

  \end{enumerate}
\end{theorem}

\begin{remark}
  \label{rmk:diskblock}
  It follows from Theorem~\ref{thm:abvGibbs} that for each
  $\gamma\in\Gamma$:
  \begin{itemize}
  \item items $(c)$ and $(d)$ ensure that
    $\leb_\gamma$-a.e. $x\in\gamma$ is $\mu$-generic by the Ergodic
    Theorem;
  % \item item (vi) ensures that $\leb_\gamma$-a.e. $x\in\gamma$ belongs
  %   to $H$ via Proposition~XXX;
  \item item $a(ii)$ implies that $\gamma$ is contained in the block
    $B_{\U}$ for the schedule $\U(x)$ of the
    $\sigma_u^{1/2}$-hyperbolic times of $x$ for $\leb$-a.e.
    $x\in H\subset U$. Hence, we may assume that $\gamma\subset H_u$
    after perhaps slightly decreasing the value of $c_u>0$.
  \end{itemize}
\end{remark}

\begin{remark}[crossing $cu$- and $cs$-disks]
  In what follows we say that a $cu$-disk \emph{crosses} $\cC_0$ if it
  intersects the cylinder $\cC_0$ and contains a graph over
  $D^u$. Analogously, we say that a $cs$-disk \emph{crosses} $\cC_0$ if it
  intersects the cylinder $\cC_0$ and contains a graph over
  $D^s$.
\end{remark}

\subsubsection{The weakly dissipative case with one-dimensional
  center-stable direction}
\label{sec:weakly-dissip-case}

Here we obtain the first part of the statement of
Corollary~\ref{mcor:abv}.

Coupling the non-uniform expansion assumption~\eqref{eq:NUE0} with the
weakly dissipative assumption, together with one-dimensional
center-stable direction, enables us to show that for any $cu$-disk
$\gamma$ in $U$ the points $H_u\cap\gamma$ are non-uniformly contracting
along the center-stable direction. This is the \emph{mostly
  contracting} property of a dominated splitting introduced by Bonatti
an Viana in~\cite{BoV00}; see also~\cite{vasquez2006} and
Theorem~\ref{thm:Vasquez}.

Indeed, the domination assumption ensures that the angle between
$E^{cs}$ and $E^{cu}$ is uniformly bounded below away from zero and so
we find a constant $0<\kappa\le1$ so that
\begin{align*}
  |\det (Df\mid E^{cs}_x)| \cdot |\det (Df\mid E^{cu}_x)|
  \le
  \kappa\cdot |\det Df(x)|,\quad x\in\Lambda.
\end{align*}
From $s=\dim E^{cs}=1$ and weak dissipativeness we obtain
\begin{align}\label{eq:weakdom}
  \|Df^n\mid E^{cs}_x\| \cdot |\det (Df^n\mid E^{cu}_x)|
  \le
  \kappa|\det Df^n(x)|
  \le
  \kappa, \quad n\ge1, x\in\Lambda.
\end{align}
For any point $x\in H_u$ satisfying~\eqref{eq:NUE0} we can write with
$d_{cu}=\dim E^{cu}\ge1$
\begin{align}
  S_k\phi^{cs}(x)
  &\le
  k\log\kappa+\log\big|\det(Df^k\mid E^{cu}_x)^{-1}\big|
  \le
  k\log\kappa  +  d_{cu}\cdot S_k\phi^{cu}(x) \nonumber
  \\
  &\le\label{eq:NUC2}
    k\big(\log\kappa+(d_{cu}/k)S_k\phi^{cu}(x)\big)
\end{align}
and since $\kappa\in(0,1]$ we obtain~\eqref{eq:NUC} for $x\in H_u$.
That is,~\eqref{eq:NUE0} implies~\eqref{eq:NUC} in the setting of
Corollary~\ref{mcor:abv}, i.e., $H_u\subset H_s$ for
$c_s=-\log\kappa + d_{cu} c_u$.

In particular, the inequality~\eqref{eq:NUC2} ensures that $x\in H_u$
admits infinitely many \emph{simultaneous hyperbolic times};
see~\cite[Proposition 6.4]{ABV00}. It follows from this that every
disk $\gamma\in\Gamma$ is such that each $y\in\gamma$ satisfies
$\|Df^{-k}\mid E^{cs}_y\|\ge e^{k c_s/2}$ for all $k\ge1$.

Therefore, the $\mu$-generic points of $y\in\gamma$ (which are also
Oseledets regular points) have a negative Lyapunov exponent along the
central-stable direction. Thus, $\mu$ is a physical/SRB probability
measure and a $cu$-Gibbs state.

This is enough to obtain the statement of existence of finitely many
ergodic physical/SRB probability measures of Corollary~\ref{mcor:abv},
following the proof of~\cite[Proposition 6.4]{ABV00}.

% In particular, from Remark~\ref{rmk:diskblock}, $\leb_{\Sigma_0}$-a.e
% $x\in\Sigma_0$ is $\mu$-generic and thus $\lambda_{cu}^->c_u/2>0$.

% Hence, from the same arguments of
% Subsection~\ref{sec:hyperb-dominat-measu-1}, there exists $N\in\ZZ^+$
% so that $g=f^N$ is non-uniformly expanding $\mu$-a.e and thus also for
% $\leb_{\Sigma_0}$-a.e. $x$, that is
% \begin{align*}
%   \limsup\nolimits_{k\nearrow\infty}S^g_k\phi^{cu}_N(x)/k<-c_u/2.
% \end{align*}
% Therefore, we obtain for $\leb_{\Sigma_0}$-a.e. $x$ as in~\eqref{eq:NUC2}
% \begin{align*}
%   \lim_{k\to\infty}
%   \log\|Df^k\mid E^{cs}_x\|^{1/k}
%   &=
%   (1/N)\lim_{k\to\infty}\log\|Dg^k\mid E^{cs}_x\|^{1/k}
%   \\
%   &\le
%   (1/N)
%   \limsup\nolimits_{k\nearrow\infty}S_k^g\phi_N^{cs}(x)/k
%   <
%   -u\cdot c_u/2N <0.
% \end{align*}
% Since such points are $\mu$-generic, we have recovered the main
% consequence of the mostly contracting assumption in the proof of
% ~\cite[Corollary C, Section 4.1.2]{vasquez2006} which
% provides the first part of the statement of Corollary~\ref{mcor:abv}.

\subsection{Construction of the stable family}
\label{sec:recurr-cylind-over}

Proceeding with the proof of Theorem~\ref{mthm:NUHypGMY}, we assume
from now on that $\mu$ is an $f$-ergodic hyperbolic dominated
$cu$-Gibbs state. From Proposition~\ref{pr:NUEfrompos} we have
nonuniform hyperbolicity $\mu$-a.e. for a power $g=f^N$, for some
$N\ge1$. Since $\mu$ decomposes as in~\eqref{eq:ergpower} with an
$f^k$-invariant and $g$-ergodic $\nu$, and $k$ a factor of $N$, then
$\nu$ is also a $cu$-Gibbs state for $g$. Indeed, besides the positive
exponents along $E^{cu}$ we have, since $\mu$ is $cu$-Gibbs, that
\begin{align*}
  h_\nu(f^k)
  &
    =h_\mu(f^k)
  =
    k\cdot h_\mu(f)=k\int J^{cu}\,d\mu=\int S_kJ^{cu}\,d\nu
  =
  \int \log|\det Df^k\mid E^{cu}|\,d\nu 
\end{align*}
and so $\nu$ satisfies
$h_\nu(g)=\frac{N}k h_\nu(f^k)=\int \log|\det Dg\mid E^{cu}|\,d\nu$.

Thus, $g$ is nonuniformly hyperbolic on the respective ergodic basin
$B(\nu)$. Hence, after peharps replacing $f$ by some power and $\mu$
by an equivalent measure, we assume without loss that
$\leb_{\gamma}$-a.e. $x\in\gamma$ is Birkhoff generic for $\mu$ and
$\gamma\in\Gamma$ with $\wh{\Gamma}\subset\supp\mu\subset A$; and
both~\eqref{eq:NUE0} and~\eqref{eq:NUC} hold on $B(\mu)$. Therefore,
we have the assumptions of Theorem~\ref{thm:abvGibbs} and the unstable
family of disks $\Gamma$ on the cylinder $\cC_0$, which is part of a
coherent block.

% In particular,~\eqref{eq:Birkhoffgeneric} also holds.

\begin{remark}[Syncronizing returns to the coherent block $\wh\Gamma$]
  \label{rmk:syncblock}
  By assumption, $\mu$-a.e. point has $f$-coherent schedules $\wh{\U}$ of
  long reverse $\sigma_s$-hyperbolic times % and $\U_u$ of
  % $\sigma_u$-hyperbolic times, both
  with positive asymptotic
  frequency, from the results of
  Section~\ref{sec:hyperb-dominat-measu}, where
  $\sigma_s:=e^{-7c_s/8}$. %and $\sigma_u:=e^{-7c_u/8}$.
  Therefore, there exist the corresponding $f$-coherent block $B^s$
  such that $\mu(B^s)>0$.

  Since $\mu$ is $f$-invariant and ergodic, then there exists
  $\ell\ge0$ so that $\wt{H}:=\wh{\Gamma}\cap f^{-\ell}B^s$ satisfies
  $\mu(\wt{H})>0$ and $\mu$-a.e. point $x$ has positive frequency of
  visits $\cH(x)\subset\ZZ^+$ to this subset. Moreover, since $\mu$ is
  $cu$-Gibbs, we can assume without loss of generality that
  $\gamma:=W^u_x(\delta_1)\in\Gamma$ with $x\in\wh{\Gamma}$ (from
  Theorem~\ref{thm:abvGibbs}) so that $\leb_\gamma(\wt{H})>0$.
\end{remark}

\begin{figure}[htpb]
  \centering \includegraphics[width=6cm]{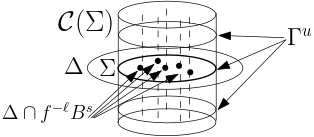}
  \caption{\label{fig:cylinder}A sketch of the center-unstable disk
    $\Delta$, the cylinder $\cC(\Sigma)$ over this disk and some
    center-stable leaves through points of $\Sigma\cap f^{-\ell}B^s$.}
\end{figure}

We write
$\cC(\Delta):=\cC_{\delta_2}(\Delta)=\cup_{x\in\Delta}W^{cs}_x(\delta_2)$
for some $0<\delta_2<\delta_1/4$ and any disk $\Delta\subset\gamma$ in
what follows; see Figure~\ref{fig:cylinder}.
We set
\begin{description}
\item[$\Sigma:=W^u_x(\delta_1/4)\subset\gamma$] a subdisk around $x$
  in the local unstable manifold through $x$ together with a small
  enough $0<\delta_2<\delta_1/4$ so that $\cC(\Sigma)\subset\cC_0$
  from Theorem~\ref{thm:abvGibbs};
\item [$\wt{H}_0:=\cC(\Sigma)\cap\wt{H}$] the subset of points of
  $\wt{H}$ inside the cylinder $\cC(\Sigma)$; see
  Figure~\ref{fig:cylinder}.  We recall that throught each
  $x\in\wt{H}$ there passes a uniformed sized stable leaf.  In
  addition
\item[$\Gamma^u:=\cC(\Sigma)\cap\Gamma$] the collection of local unstable
  manifolds of $\Gamma$ restricted to $\cC(\Sigma)$ which cross
  $\cC(\Sigma)$, and so are graphs of $C^1$ maps $\Sigma\to E^{cs}_x$
  in the local exponential chart.
\end{description}

We assume, without loss of generality, that $\mu(\wt{H}_0)>0$.

\begin{proposition}[stable lamination crosses $\cC(\Sigma)$
  $\leb-\bmod0$]
  \label{pr:fullWs}
  There exists a full $\leb_\Sigma$-measure subset $Y$ of $\Sigma$
  whose center-stable leaves $\{W^{cs}_y: y\in Y\}$ cross
  $\cC(\Sigma)$ and are uniformly contracted at a rate
  $\sigma_s^{1/2}$. Moreover, each return time $n\in\ZZ^+$ of $x\in Y$
  to $\wt{H}_0$ is a $\sigma_u^{1/4}$-hyperbolic time in the
  center-unstable direction.
\end{proposition}

\begin{remark}[stable tail condition $\&$ full disks inside
  $\wt{H_0}$]\label{rmk:fulldiskH}
We do not need to assume any tail condition on the speed of
convergence of non-uniform contraction along the center-stable
direction to obtain exponential mixing, since we obtain
\emph{uniformly long stable leaves with uniform contraction almost
  everywhere inside certain cylinders on the ambient space}. This
provides the the ``generalized horseshoe with infinitely many returns
in variable times'' which is known since Young~\cite{Yo98} to control
the speed of mixing.

Proposition~\ref{pr:fullWs} in particular ensures that $\wt{H}_0$
contains a full $\leb_\gamma$-measure subset of each
$\gamma\in\Gamma^u$ (perhaps considering smaller values of
$c_u,c_s>0$).
\end{remark}

\begin{proof}[Proof of Proposition~\ref{pr:fullWs}]
  We consider the \emph{induced transformation
    $F:\wh{\Gamma^u}\to\wh{\Gamma^u}$ given by the first return map to
    $\wh{\Gamma^u}:=\cup_{\gamma\in\Gamma^u}\gamma$ with induced time
    $\tau:\wh{\Gamma^u}\to\ZZ^+$}, that is, $F(x):=f^{\tau(x)}(x)$. This
  return map $F$ is well defined, since
  $\mu(\wh{\Gamma^u})\ge\mu(\wt{H}_0)>0$, and also bimeasurable and
  invertible, since $f$ is a diffeomorphism.

  We also consider the interated return map
  $F^i:\wh{\Gamma^u}\to\wh{\Gamma^u}, i\ge1$ and the corresponding induced
  iterated return time $\tau^i:\wh{\Gamma^u}\to\ZZ^+$ so that
  $F^i(x)=f^{\tau^i(x)}(x)$ for $\mu$-a.e. $x\in\wh{\Gamma^u}$.

  It follows, from Remark~\ref{rmk:diskblock}, that each $\tau(x)$ is
  a $\sigma_u^{1/2}$-hyperbolic time for $\mu$-a.e. $x\in\wh{\Gamma^u}$.
  Hence, from Lemma~\ref{le:contraction}, there exists a pre-disk
  $V_{\tau(x)}(x)\subset W^u_x(\delta_1)$. We consider the pre-disk
  restricted to $\wh{\Gamma^u}$, given by
  \begin{align*}
    \wt{V}_{\tau(x)}(x):=
    \big(f^{\tau(x)}\mid V_{\tau(x)}(x)\cap\wh{\Gamma^u})^{-1}
    \big(W^u_{Fx}(\delta_1)\cap\wh{\Gamma^u}\big).
  \end{align*}
  We are now ready to consider the following subset
  \begin{align*}
    Y:=\bigcup_{i\ge1} \bigcup_{x\in\wt{H}_0}
    %\left\{
    \big( f^{\tau^i(x)}\mid \wt{V}_{\tau^i(x)}(x) \big)^{-1}(\wt{H}_0).
    %\text{ for each pair }
    %x\in\wt{H}_0
    % \text{ so that } F^ix\in\wt{H}_0
    %\right\}.
  \end{align*}
  In what follows we show that: (i) $Y$ is Borel measurable; (ii) each
  of its points have long stable manifolds; (iii) each
  $\gamma\in\Gamma^u$ intersects $Y$ in a full $\leb_\gamma$-measure
  subset.

  \begin{remark}[image of pre-disks in $Y$ crosses $\Gamma^u$]
    \label{rmk:openpredisk}
    Given $x\in\wt{H}_0$ there exists $m=m(x)\in\ZZ^+$ so that if
    $\tau(x)>m$, then
    $f^{\tau(x)}\big( \wt{V}_{\tau(x)}(x) \big) =
    W^u_{Fx}(\delta_1)\cap\wh{\Gamma^u}$, that is, \emph{the image of the
      local pre-disk crosses $\Gamma^u$}, since
    $V_{\tau(x)}(x)\subset\wh{\Gamma^u}$ due to shrinking diameter when
    $\tau(x)$ grows.
  \end{remark}

  \begin{lemma}
    \label{le:measY}
    The subset $Y\subset\wh{\Gamma^u}$ is Borel measurable.
  \end{lemma}

\begin{proof}
  If we set $\wt{H}_n:=\wt{H}_0\cap f^{-n}\wt{H}_0$, then we have
  \begin{align}\label{eq:Yalt}
    Y = \bigcup_{n\ge1} \bigcup_{x\in \wt{H}_n }
     \big(  f^n\mid \wt{V}_n(x) \big)^{-1}(\wt{H}_0).
  \end{align}
  In addition, $\wt{H}_n$ is covered by the $(n+1,\delta_1)$-dynamical
  balls $\{B(x,n+1,\delta_1): x\in\wt{H}_0\}$, from
  Remark~\ref{rmk:prediskdynball}. Then, since the ambient space $M$
  is a smooth manifold, we can find a \emph{denumerable subcover}
  $\{B(x_k,n+1,\delta_1):k\ge1\}$ where $(x_k)_{k\ge1}$ is a sequence
  in $\wt{H}_0$. Again from Remark~\ref{rmk:prediskdynball} the family
  $$\{B(x_k,n+1,\delta_1)\cap\cC(\Sigma): k\ge1\}$$ covers all the
  pre-balls $V_n(x)$ contained in the countable
  union~\eqref{eq:Yalt}. Thus
  \begin{align*}
    \bigcup_{x\in \wt{H}_n} \big(  f^n\mid \wt{V}_n(x) \big)^{-1}(\wt{H}_0)
    =
    \bigcup_{k\ge1} \big(f^n\mid B(x_k,n+1,\delta_1) \big)^{-1}(\wt{H}_0)
  \end{align*}
  and $Y$ can be rewritten as a countable union of clearly
  Borel measurable subsets. Hence, $Y$ is Borel measurable.
\end{proof}

\begin{lemma}[Uniformly long stable leaves through $Y$]
  \label{le:inducehyptimes}
  Let $i\ge 1$ be a return time to $\wt{H}_0$ of $x\in\wt{H}_0$ and
  $\wt{V}:=\wt{V}_{\tau^i(x)}(x)$ the hyperbolic pre-disk in
  $\Gamma^u$.  Then
  \begin{enumerate}
  \item whenever $n=\tau^i(x)>\ell$, then $n$ is a reverse
    $\sigma_s^{1/2}$-hyperbolic time of each $y\in\wt{V}$ with
    respect to $\ell$, that is,
    $\| Df^k\mid E^{cs}_{f^\ell y}\|\le (\sigma_s^{1/2})^k, 1\le k\le
    n-\ell$;
  \item we can find $0<\delta_2<\delta_1/4$ small enough so that, for
    each $y\in\wt{V}$, the center-stable leaf $W^{cs}_y(\delta_2)$
    satisfies item (2) of Proposition~\ref{pr:reverhyp} with
    contraction rate $\sigma=\sigma_s^{1/2}$ modulo a uniform
    constant.
  \end{enumerate}
  More precisely, there exists $\delta_2>0$ so that each $y\in Y$
  admits a stable leaf $W^s_y(\delta_2)$ of uniform size, tangent to
  the center-stable subbundle and with uniform rates of contraction:
  there exists $C_s>0$ such that for for all $z\in W^s_y(\delta_2)$ we
  have
  \begin{align*}
    \dist(f^kz,f^ky)\le C_s \big(\sigma_s^{1/4}\big)^k \dist(z,y),
    \quad\forall k\ge0.
  \end{align*}
  In addition, we also have that
  \begin{enumerate}[resume]
  \item $n$ is a $\sigma_u^{1/4}$-hyperbolic time for each
    $z\in W^s_y(\delta_2)$ and any $y\in\wt{V}$.
  \end{enumerate}
\end{lemma}

% Figura ilustrativa aqui!!

\begin{figure}[htpb]
  \centering
  \includegraphics[width=10cm]{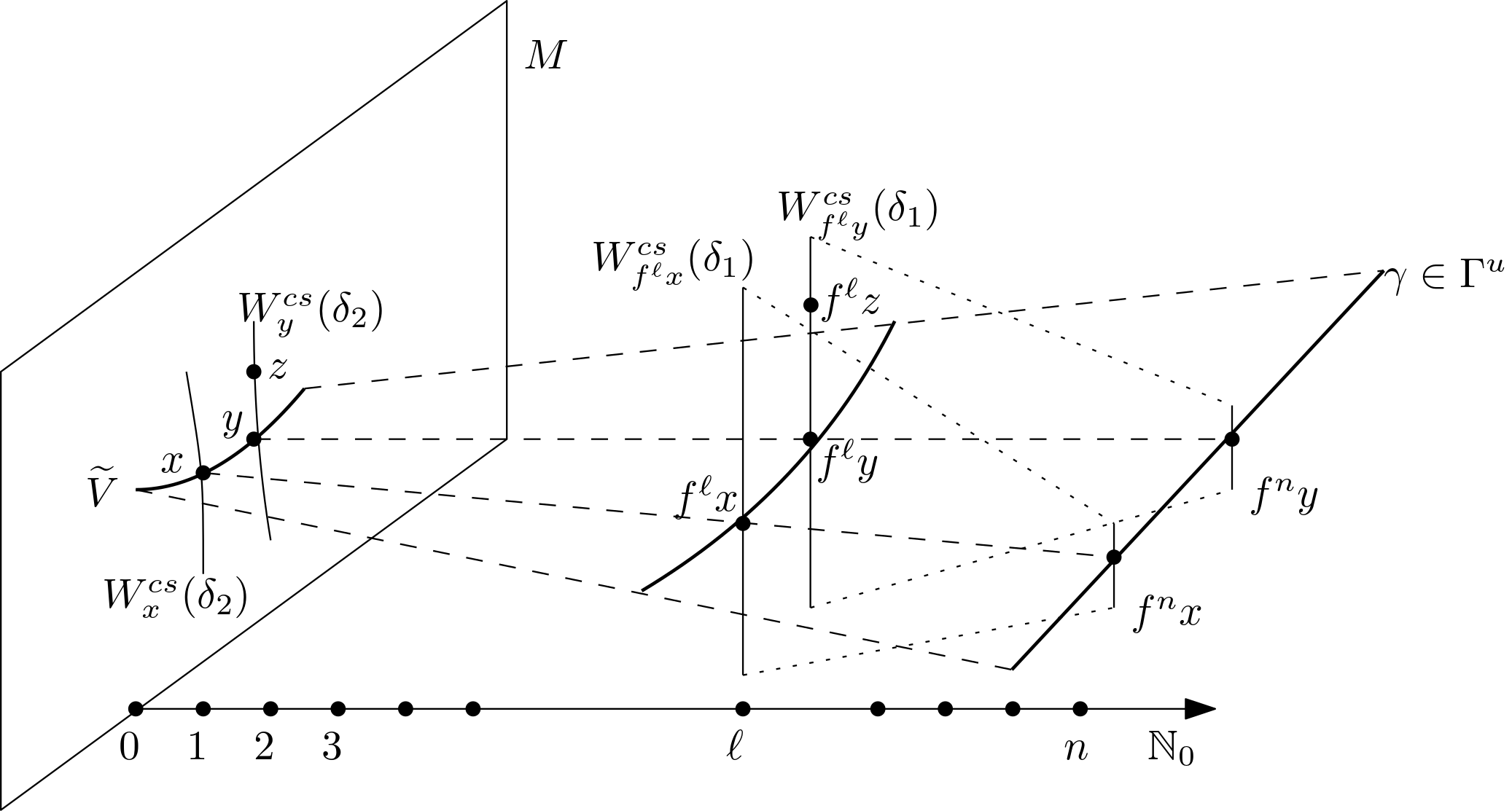}
  \caption{\label{fig:uniflongstbl} Construction of
    long stable leaves at $y\in\wt{V}$.}
\end{figure}

\begin{proof}
  For item (1), we note that $\wt{V}$ contains $x\in\wt{H}_0$ so
  $\ell$ is a long $\sigma_s$-hyperbolic time for $x$, and
  \begin{itemize}
  \item $S_k\phi^{cs}(f^\ell x)<-7c_s k/8$ for all $k\ge1$ since
    $\sigma_s=e^{-7c_s/8}$; and
  \item $\dist_{f^i\omega}(f^iy,f^ix)\le\delta_1/4$ for $y\in\wt{V}$
    and $i=0,\dots,n$, by definition of hyperbolic pre-disk $\wt{V}$
    contained in $\Gamma^u$, since $n$ is a
    $\sigma_u^{1/2}$-hyperbolic time for all $y\in \wt{V}$.
  \end{itemize}
  Then, by the choice of $\delta_1$ in~\eqref{eq:delta1}, together
  with Remark~\ref{rmk:delta1simult}, we get
  \begin{align*}
    S_k\phi^{cs}(y)
    \le
    S_k\phi^{cs}(x)
    +
    7 k c_s/16;
    \quad k=1,\dots, n; \quad y\in\wt{V}.
  \end{align*}
  Hence, if $n>\ell$, then
  $ S_k\phi^{cs}(f^\ell y) \le S_k\phi^{cs}(f^\ell x) + (7 c_s/16)k \le
  (-7 c_s/16)k$, $k=1,\dots,n-\ell$, and this shows that each
  $y\in\wt{V}$ has $\ell$ as a reverse $\sigma_s^{1/2}$-hyperbolic
  time with respect to $n$, as stated.

  For item (2), from Proposition~\ref{pr:longstable}, the $cs$-disk
  $W^{cs}_{f^\ell y} (\delta_1)$ through $f^\ell y$ is uniformly
  contracted during the next $n-\ell$ iterates at the rate
  $\sigma_s^{1/4}$; see Figure~\ref{fig:uniflongstbl}. Since $f$ is
  $L$-Lipschitz, we can find $0<\delta_2<\delta_1/4$ so that
  $f^\ell\big( W^{cs}_y(\delta_2)\big) \subset W^s_{f^\ell
    y}(\delta_1)$. Then, for each $z\in W^{cs}_y(\delta_2)$ and
  $k>\ell$
  \begin{align*}
    \dist( f^kz , f^k y)
    &\le
      \dist_{f^kW^{cs}_y(\delta_2)}\big( f^{k-\ell}f^\ell z , f^{k-\ell}
      f^\ell y \big)
      \le
      (\sigma_s^{1/4})^{k-\ell}\dist_{f^\ell W^{cs}_y(\delta_2)}(f^\ell
      z,f^\ell y)
    \\
    &\le
      L^\ell (\sigma_s^{1/4})^{k-\ell}\dist_{W^{cs}_y(\delta_2)}(z,y)
      \le
      C_1\big(L\sigma_s^{-1/4}\big)^\ell (\sigma_s^{1/4})^k \dist(y,z)
  \end{align*}
  where, in the last inequality, we used the bound on curvature of all
  $cs$- and $cu$-disks; see
  Subsection~\ref{sec:conseq-existence-dom}. If $0<k\le\ell$, then
  \begin{align*}
    \dist(f^kz,f^ky)
    &\le
      L^k\dist(y,z)
      =
      \big(L\sigma_s^{-1/4}\big)^k\cdot (\sigma_s^{1/4})^k \dist(z,y).
  \end{align*}
  If we set
  $C_s:=\max\{1, C_1(L/\sigma_s^{1/4})^\ell, (L/\sigma_s^{1/4})^i:
  i=1,\dots,\ell-1\}$, then we deduce the bound stated in item (2) for
  all $k\ge0$.

  For item (3), for $z\in W^s_y(\delta_2)$ and $y\in\wt{V}$, we have
  the following for each $k>0$
  \begin{align*}
    \dist_{f^kW^s_y(\delta_1)}(f^ky,f^kz)\le
    C_s\delta_2(\sigma_s^{1/4})^k\le C_2\delta_2 \le \delta_1
  \end{align*}
  if we let $\delta_2\le \delta_1/C_2$, from item (2); see
  Figure~\ref{fig:uniflongstbl}.  Moreover, from item (1) and the
  choice~\eqref{eq:delta1}, together with
  Remark~\ref{rmk:delta1simult}, we have
  $\phi^{cu}(f^ky)-\phi^{cu}(f^kx) < \log\sigma_u^{-1/4}= 7c_u/32$ and
  so
  \begin{align*}
    S_{n-k}\phi^{cu}(f^kz)
    <
    S_{n-k}\phi^{cu}(f^ky)
    + 7(n-k)c_u/32
    \le
    -7(n-k)c_u/32,
    \quad 0\le k<n.
  \end{align*}
  Thus, $n$ becomes a $\sigma_u^{1/4}$-hyperbolic time for $z$.  The
  proof is complete.
\end{proof}

The following result ensures that $Y$ has full measure inside
$\cC(\Sigma)$.

\begin{lemma}
  \label{le:Finv}
  The subset $Y$ is forward $F$-invariant $F(Y)\subset Y$ with
  positive $\mu$-measure.
\end{lemma}

Since $f$ is invertible and bimeasurable, then $F$ is also
invertible. Moreover, since $\mu$ is $f$-invariant and ergodic, then
the normalized restriction $\mu_0$ of $\mu$ to $\wh{\Gamma^u}$ is
$F$-invariant and ergodic. Therefore, from Lemma~\ref{le:Finv}, we
conclude that $\mu(\wh{\Gamma^u}\setminus Y)=0$, that is,
$Y=\wh{\Gamma^u}, \mu\bmod 0$. Hence, considering the absolutely
continuous disintegration of $\mu$ along the leaves of $\Gamma^u$
(cf. item (d) of Theorem~\ref{thm:abvGibbs}) we obtain
$\leb_\gamma(\wh{\Gamma^u}\setminus
Y)=\mu_\gamma(\wh{\Gamma^u}\setminus Y)=0$ for $\hat\mu$-almost all
$\gamma\in\Gamma^u$.

Thus, we can assume, without loss of generality, that
$\leb_{\Sigma}$-a.e. point $x\in\Sigma$ admits a uniformly sized
stable leaf $W^s_x(\delta_2)$ with uniform rate of forward
contraction. This completes the proof of Proposition~\ref{pr:fullWs}
assuming Lemma~\ref{le:Finv}.
\end{proof}

We are left to provide the following.

\begin{proof}[Proof of Lemma~\ref{le:Finv}]
  Clearly $Y\subset\wh{\Gamma^u}$ and $x\in\wt{H}_0$ returns to $\wt{H}_0$
  in some iterate $k>0$, thus we obtain  $x\in(f^{\tau^k(x)}\mid
  \wt{V}_{\tau^k(x)}(x))^{-1}(\wt{H}_0)$. Hence, $Y\supset\wt{H}_0$ and it
  follows that $\mu(Y)>0$, by construction of $\wt{H}_0$. We are left
  to prove the forward $F$-invariance of $Y$.

  Let $y\in Y$ be given. Then there exist $x\in \wt{H}_0$ and $k\ge1$
  so that $f^{\tau^k(x)}(y)\in\wt{H}_0$. Hence, there exists
  $\ell\ge1$ so that $F^\ell(y)\in\wt{H}_0$.

  If $\ell=1$, then $Fy\in\wt{H}_0\subset Y$. Otherwise, we have
  $\ell>1$ and $F^{\ell-1}(Fy)\in\wt{H}_0$. Hence, by definition of
  $F$ we have
  \begin{itemize}
  \item $Fy\in\wh{\Gamma^u}$ and
    $Fy\in\wh{V}_{\tau^{\ell-1}(Fy)}(Fy)$; and also
    $F^{\ell-1}(Fy)\in\wt{H}_0$; and moreover
  \item $\tau^{\ell-1}(Fy)$ is a $\sigma_u^{1/2}$-hyperbolic time for $Fy$, by
    definition of $\wt{H}_0$.
  \end{itemize}
  Thus, by definition of $Y$, we conclude that $Fy\in Y$, completing
  the proof.
\end{proof}

% \hrulefill

% \begin{remark}[Absolutely continuous holonomy]
%   \label{rmk:abscontholo}
%   This construction provides a stable holonomy
%   $\pi:\cC_0\to\Gamma_0$ which is absolutely continuous with bounded
%   density; see the last argument in the proof of
%   Theorem~\ref{mthm:PesinC1+}.
% \end{remark}

% Since the stable leaves
% through $\wh{H}_0$ are a restriction of the continuous lamination of
% the center-stable leaves forming $\cC(\Delta)$, we also write
% $\pi:\cC(\Sigma)\to\Sigma$ for the projection on $\Sigma$ given by
% the holonomy along center-stable leaves.

% \begin{remark}[H\"older continuous holonomy]
%   \label{rmk:Holderholo}\margem{New subsection of Sec.4?}
%   Using the ``H\"older invariant section theorem'' from \cite[Theorem
%   5.18]{Sh87} and \cite[Section 3.2]{Wi98} we have that the locally
%   defined center-stable lamination on $\cC(\Sigma)$ induces
%   $\alpha$-H\"older holonomies along its leaves between $cu$-disks
%   contained in $\cC(\Sigma)$ for some $0<\alpha<1$, which are a
%   H\"older-continuous extension of the stable holonomy from
%   Remark~\ref{rmk:abscontholo}. 
% \end{remark}

\subsection{The full GMY structure with integrable
  return times}
\label{sec:existence-full-gmy}

We are now ready to present the following.

\begin{proof}[Proof of Theorem~\ref{mthm:NUHypGMY}]
  We have already defined the family $\Gamma^u$ of unstable manifolds
  and set $Y_0:=\Sigma\cap Y$, a full $\leb_\Sigma$-measure subset of
  the local unstable manifold $\Sigma\in\Gamma^u$, where $Y$ is given by
  Proposition~\ref{pr:fullWs}, and the family
  $\Gamma^s:=\{W^s_x(\delta_2): x\in\Sigma\}$ of local stable
  manifolds. Both these families are a subset of the respective
  families of center-unstable and center-stable manifolds given by the
  dominated splitting and, thus, $\Gamma^u$ and $\Gamma^s$ are
  automatically continuous.

  We show that $\Lambda:=(\cup\Gamma^u)\cap(\cup\Gamma^s)$ has a
  hyperbolic product structure with respect to an induced return map
  under $f$. By the previous contructions we already have conditions
  (1)-(4) of Subsection~\ref{sec:gibbsm-struct} from the definition of
  GMY structure, together with item (I) for $\gamma=\Sigma$.

  In order to define the Markov return map, we consider the sequence
  of subsets
  $H_{-n}:=f^{-n}(\wt{H}_0)$ for each $n\ge1$. % We recall that $m$ as
  % a $\sigma_u^{1/2}$-hyperbolic time and admits the corresponding
  % pre-ball $V_m(f^{n-m}x)$ for each $x\in \wt{H}_n$ and all $m>0$.
  There exists $\wt{\theta}\in(0,1]$ so that for $\leb_\Sigma$-a.e. $x$
  we can find $n_0\in\ZZ^+$ satisfying
 \begin{align*}
   n\ge n_0 \implies
   \#\big\{1\le j\le n : x\in H_{-j}\big\}
   =
   \#\big\{1\le j\le n : f^jx\in \wt{H}_0\big\}
   >
   n\wt{\theta},
 \end{align*}
 and so we can define
 \begin{align}\label{eq:htilde}
   \wt{h}_\theta(x):=\min\left\{N\ge1:
   \#\big\{1\le j\le n : x\in H_{-j} \big\}\ge n\wt{\theta},
   \forall n\ge N\right\},
 \end{align}
 where $\wt{\theta}=\wt{\theta}(\sigma_u^{1/2})$, given by
 Lemma~\ref{le:pliss} of Pliss, depends only on $f$ and on the rate
 $\sigma_u^{1/2}=e^{-7c_u/16}$.  We are ready to obtain the following.

  \begin{theorem}{\cite[Proposition 7.16 \;\&\; Theorem
      5.1]{Alves2020b}}\label{thm:tailGMY}
    Given $N_0\ge1$ there exists a $\leb_\Sigma-\bmod0$ partition $\cP$
    of $\Sigma$ into domains $\omega_n$ so that
    $\omega_n\subset V_n(x)$ for some $x\in H_{-n}$ and $n\ge N_0$.
    Setting $R(x)=n$ for $x\in\omega_n\in\cP$, we get that
    \begin{enumerate}
    \item for every $n\ge1$ there are finitely many $\omega\in\cP$
      with $R(\omega)=n$;
    \item $f^R\mid\omega:\omega\to W^u_{f^Rx}(\delta_1)\cap\wh{\Gamma^u}$
      maps each $\omega\in\cP$ to an unstable leaf crossing $\cC(\Sigma)$;

    \item there are $(S_i)_{i\ge1}$ subsets of $\Sigma$ so that
      $\sum_{n\ge1}\leb_\Sigma(S_n)<\infty$ and $H_{-n}\cap\{R>n\}\subset
      S_n$ for all $n\ge1$.
    \item there are $(E_i)_{i\ge1}$ subsets of $\Sigma$ so that
      $\leb_\Sigma(E_i)$ tends to zero exponentially fast and
      $\{R>n\}\subset\{\wt{h}_\theta>n\}\cap E_n$ for all $n\ge1$.
    \end{enumerate}
  \end{theorem}

  \begin{proof}
    This is essentially the statement of~\cite[Proposition
    7.16]{Alves2020b}.  Since, in our setting, we already have an
    ergodic physical/SRB measure, we know that
    $\leb_\Sigma$-a.e. $x\in\Sigma$ is $\mu$-generic. Thus,
    $\leb_\Sigma$-a.e. $x$ belongs to infinitely many subsets from
    $(H_{-n})_{n\ge1}$.  The full statement of~\cite[Proposition 7.16
    \;\&\; Theorem 5.1]{Alves2020b} demands an extra $(I_3)$ condition
    and \cite[Lemma 7.15]{Alves2020b}, which out setting automatically
    provides with the constants $L=\ell=0$, in the notation
    of~\cite[Chapters 5 \;\&\; 7]{Alves2020b}.
  \end{proof}

  We set $R\mid \Lambda_i \equiv R_i = R(\omega_i)$ for $i\ge1$, where
  \begin{align*}
    \Lambda_i:=\Gamma^u\cap  \bigcup_{x\in \omega_i \cap Y_0  }
    W^s_x(\delta_2).
  \end{align*}
  Then, to obtain that
  $\Lambda=\cup_i\Lambda_i$ has full GMY structure with recurrence
  time $R$, we follow verbatim the proof of~\cite[Proposition
  7.21]{Alves2020b}, since in our setting we have
  \begin{enumerate}[(i)]
  \item the function
    $x\in\Lambda_i\mapsto\log|\det Df\mid T_{f^kx}f^k\gamma|$ is
    $(L_1,\zeta)$-H\"older-continuous for all $0\le k<R_i$, from
    Proposition~\ref{pr:curvature} and Corollary~\ref{cor:distortion};
  \item uniform contraction of the stable leaves from $\Gamma^s$
    covering the cylinder $\cC(\Sigma)$ from
    Proposition~\ref{pr:fullWs};
  \item the subbundles $E^{cs}$ and $E^{cs}$ are H\"older-continuous,
    from the domination assumption.
  \end{enumerate}
  Following the arguments in~\cite[Proposition 7.21]{Alves2020b} we
  obtain all the conditions (I)-(VI) with each disk of $\Gamma^u$
  contained in $\Lambda$.

  The integrability of the recurrence time $R$ follows from the
  arguments in~\cite[Section 7.3 of Chapter 7]{Alves2020b}.  This
  completes the argument for the existence of the GMY structure with
  integrable return time and finishes the proof of
  Theorem~\ref{mthm:NUHypGMY}.
\end{proof}

At this point we are able to complete the following.

\begin{proof}[Proof of Corollary~\ref{mcor:abv}]
  From the first part of the statement of Corollary~\ref{mcor:abv},
  obtained in Subsection~\ref{sec:weakly-dissip-case}, we have
  finitely many $\mu_1,\dots,\mu_k$ ergodic physical/SRB measures
  which are $cu$-Gibbs states. Hence we are in the setting of
  Theorem~\ref{thm:abvGibbs} for each $\mu_i$ and the second part of
  the statement of Corollary~\ref{mcor:abv} follows. For the equality
  between geometric, ergodic and topological basins, see the proof of
  Theorem~\ref{mthm:basin} in the next
  Section~\ref{sec:geometr-basin-coinci}.
\end{proof}

%%%%%%%%%%%%%%%%%%%%%%%%%%%%%%%%%%%%%%%%%%%%%%%%%%%%%%%%%%%%

\section{Speed of mixing from the GMY structure}
\label{sec:speed-mixing-from}

To prove Theorem~\ref{mthm:Gibbsmix} we recall the following standard
result.

\begin{theorem} \cite[Theorem 4.15]{Alves2020b}
  \label{thm:corrdecay} Let $f : M \circlearrowleft $ be a $C^{1+\eta}$
  diffeomorphism, for some $0<\eta\le1$, admitting a GMY structure
  $\Lambda$ with integrable recurrence time $R:\Lambda\to\ZZ^+$ and
  $\mu$ be the unique ergodic physical/SRB measure for $f$ with
  $\mu(\Lambda) > 0$. 
  If $\gcd(R)=q$, then $f^q$ has $p\le q$ exact invariant probability
  measures $\mu_i, i=1,\dots,p$ so that $f_*\mu_i=\mu_{(i+1)\mod p}$
  and $p\cdot \mu=\sum_{i=1}^p\mu_i$. Moreover, for all such $i$ and $n>1$
  \begin{enumerate}
  \item if $\leb_\gamma\{R\ge n\}\le Cn^{-\alpha}$ for some 
    $\gamma\in\Gamma^u$, $C>0$ and $\alpha>1$, then for all
    $\eta$-H\"older observables $\vfi,\psi:M\to\RR$ there is $C'>0$ so
    that
    $\corr_{\mu_i}(\vfi,\psi\circ f^{qn})\le C'n^{-\alpha+1}$.
  \item if $\leb_\gamma\{R\ge n\}\le C e^{-c n^\alpha}$ for some
    $\gamma\in\Gamma^u$, $C,c>0$ and $0<\alpha\le1$, then there exists
    $c'>0$ so that for $\eta$-H\"older observables $\vfi,\psi:M\to\RR$
    there is $C'>0$ for which
    $\corr_{\mu_i}(\vfi,\psi\circ f^{qn})\le C'e^{ - c' n^\alpha}$.
  \end{enumerate}
\end{theorem}

We relate the tail of return times $R$ with the expansion time
function $h$ to obtain the following.

\begin{proof}[Proof of Theorem~\ref{mthm:Gibbsmix}]
  The first statement of Theorem~\ref{mthm:Gibbsmix} is a consequence
  of Theorem~\ref{mthm:PesinC1+}, providing the power $g=f^N$ with a a
  physical/SRB measure for $g$. Then Theorem~\ref{mthm:NUHypGMY}
  ensures the existence of a GMY structure.

  Let us fix $\gamma\in\Gamma^u$ contained in GMY structure.  We claim
  that condition (1) or (2), of the statement of
  Theorem~\ref{thm:corrdecay}, holds whenever the tail condition on
  $h$ stated in items (1) and (2) of Theorem~\ref{mthm:Gibbsmix}
  holds, respectively.

  To prove the claim, we recall the definition of the tail function
  $h(x)$ from~\eqref{eq:cutail} and consider
  \begin{align*}
    h_\theta(x)
    :=
    \min\left\{N\ge1: \#\{1\le i\le n: x\in H_i\}\ge n \theta_1, 
    \forall\,n\geq N\right\},
  \end{align*}
  where we write $H_i=\{x\in M: i$ is a $\sigma_u^{3/7}$-hyperbolic
  time for $x \}$ and
  \begin{align*}
    \theta_1
    :=
    \frac{e^{c_u/2}-e^{3c_u/8}}{\sup(-\phi^{cu})-e^{3c_u/8}}
    <
    \theta_0
    :=
    \frac{\log\sigma_u^{-1/2}-\log\sigma_u^{-3/7}}{L^u-\log\sigma_u^{3/7}}
    =
    \frac{e^{7c_u/8}-e^{3c_u/8}}{\sup(-\phi^{cu})-e^{3c_u/8}},
  \end{align*}
  a lower bound for the frequency provided by Lemma~\ref{le:pliss} of
  Pliss.

  \begin{remark}
    \label{rmk:wtH}
    Note the subtle difference between $H_i$ and $H_{-i}$ from
    Subsection~\ref{sec:existence-full-gmy}, and also between
    $\wt{h}_\theta$ from~\eqref{eq:htilde} and $h_\theta$, in what
    follows.
  \end{remark}

  We recall, from the proof of Theorem~\ref{mthm:NUHypGMY}, that
  $R(x)=n$ means that $f^n(x)\in\wt{H}_0$, and so $n$ is a
  $\sigma_u^{1/2}=(e^{-7c_u/16})$-hyperbolic time for $x$. In
  particular, we have
  \begin{align*}
    S_n\phi^{cu}(x) < -7c_u n/8 = n \log \sigma_u^{1/2}.
  \end{align*}
  Using Lemma~\ref{le:pliss} of Pliss with the rates
  $c_2=\log\sigma_u^{-1/2}>c_1=\log\sigma_u^{-3/7}$, there are
  $\ell \ge \theta_0n$ iterates $1\le n_1<\dots<n_\ell <n$ which are
  $\sigma_u^{3/7}(=e^{-3c_u/8})$-hyperbolic times for $x$.

  This means that each visit to $\wt{H}_0$ at time $n$ ensures the
  existence of at least $\theta_0n$ previous
  $\sigma_u^{1/3}$-hyperbolic times. Thus, we get
  \begin{align}\label{eq:tailA}
    \{\wt{h}_\theta > n \} \subset \{ h_\theta > n\}.
  \end{align}
  Moreover, if $h(x)=N$, then $S_n\phi^{cu}(x)<-c_un/2$ for all
  $n\ge N$ by definition of $h(x)$ in~\eqref{eq:cutail}. Using again
  Lemma~\ref{le:pliss}, we can find $\ell\ge \theta_1 n$ times
  $1\le n_1<\dots<n_\ell\le n$ which are $\sigma_u^{3/7}$-hyperbolic
  times for $x$.

  This shows that $h_\theta(x)\le N$, since the $\theta_1$-frequency of
  $\sigma_u^{3/7}$-hyperbolic times is achieved at least from time $N$
  onwards. Hence, we arrive at
  \begin{align}\label{eq:tailB}
    \{h_\theta>N\}\subset \{h>N\}.
  \end{align}
  Altogether, from~\eqref{eq:tailA} and~\eqref{eq:tailB}, we obtain
  $\{\wt{h}_\theta>N\}\subset\{h>N\}$. Thus, the tail of $R$ in
  Theorem~\ref{thm:tailGMY} is given by the tail of $\wt{h}_\theta$
  which, in turn, is given by the tail of $h$.

  From item (4) of Theorem~\ref{thm:tailGMY}, the tail of $R$
  satisfies the conditions of items (1) or (2) of
  Theorem~\ref{thm:corrdecay} (that is, polinomial or (sub)exponential
  decay) if the tail of $h$ satisfies the same
  conditions.
  The proof of Theorem~\ref{mthm:Gibbsmix} is complete.
\end{proof}

%%%%%%%%%%%%%%%%%%%%%%%%%%%%%%%%%%%%%%%%%%%%%%%%%%%%%%%%%%%%%

\section{Geometric and ergodic basins coincide}
\label{sec:geometr-basin-coinci}

Here we prove Theorem~\ref{mthm:basin} as a corollary of
Theorem~\ref{mthm:NUHypGMY}. Since the statement of these results have
the same assumptions, we can assume that we have a GMY structure for
the ergodic hyperbolic dominated $cu$-Gibbs state $\mu$.

This is given by a cylinder $\cC(\Sigma)$ with positive measure, over
an unstable disk $\Sigma\subset A$, such that $\leb_\Sigma$-a.e. $x$
is $\mu$-generic and the corresponding stable manifold $W^s_x$
contains a stable leaf crossing the cylinder.  Moreover, the family of
stable leaves $W^s(\Sigma)=\{ W^s_x : x\in\Sigma\}$ contains a full
volume subset $W$ of $\cC$. Each element of $W$ is positively
assymptotic with the positive trajectory of some $\mu$-generic point
of $\Sigma\in A$. Thus, $W$ is contained in the geometric basin $G(A)$
of the attracting set $A$, by construction, and also in the ergodic
basin $B(\mu)$ and topological basin $B(A)$.

Let $B=B(p,\delta)$ be a ball in the interior of $\cC$ and
$\vfi:M\to[0,+\infty)$ be a non-negative continuous observable
supported in $B$ with $\mu(\vfi)>0$. Then, for any $y\in B(\mu)$
we have
\begin{align*}
  \wt{\vfi}(y)=\lim_{n\to+\infty}\frac1nS_n\vfi(y)=\mu(\vfi)>0
\end{align*}
and so there exists $n\in\ZZ^+$ so that
$f^ny\in B\subset\interior(B)$. Therefore, we can find a neigborhood
$V$ of $y$ so that $f^nV\subset B$ and so, because $f$ is a
diffeomorphism, the preimage of $W$ fills a full volume subset of $V$:
$\leb\big ( V\setminus (V\cap f^{-1}W) \big) = 0$.

This shows that a neighborhood of any point of the ergodic basin
contains a full volume subset of simultaneously the geometric basin,
ergodic basin and topological basin. We conclude that these basins
coincide over the ergodic basin $B(\mu)$ of $\mu$, that is
\begin{align*}
  B(\mu)=G(\supp\mu) \big( \subset B(A) \big)
\end{align*}
except, perhaps, a zero volume susbset of points.

In case $\leb(U\setminus H)=0$, since $B(A)\supset U\cup G(A)$ (by
definition of attracting set) and $G(A)\supset G(\supp\mu)$ for each
ergodic hyperbolic $cu$-Gibbs state, then from
Theorem~\ref{thm:Vasquez} we deduce
$B(A)=B(\mu_1)\cup\ldots\cup B(\mu_k)=G(\supp\mu_1)\cup\ldots\cup
G(\supp\mu_k)\subset G(A)\subset B(A)$ and so we have equality
throughout (perhaps except a zero Lebesgue measure subset).  This
completes the proof of Theorem~\ref{mthm:basin} and the basin claim of
Corollary~\ref{mcor:abv}.

%%%%%%%%%%%%%%%%%%%%%%%%%%%%%%%%%%%%%%%%%%%%%%%%%%%%%%%%%%%%%

\def\cprime{$'$}

\bibliographystyle{abbrv}
\bibliography{../bibliobase/bibliography}

\end{document}